%% file: sigproc09.tex
\documentclass[onecolumn]{elsart3p}
\usepackage[T1]{fontenc}
\usepackage{amsmath, epsfig,amssymb, graphicx}
\usepackage{setspace}
\usepackage{algorithmic, appendix, enumitem}
\graphicspath{{figures/}}

\input jidef
\input jidef_en

\input alphabet         
\input beginend   
\def\MT{^{-\tD}}
\usepackage{color}
\definecolor{col1}{rgb}{.7,0,0}
\def\suppJI#1{{\footnotesize \setstcolor{blue}\st{#1}}}
\def\remJI#1{{\noindent\color{blue}{{\footnotesize JI: #1}}}}
\def\addJI#1{{\noindent\color{blue}{#1}}}
\renewcommand{\algorithmiccomment}[1]{\!\texttt{\small{\% #1}}}

\newcommand{\mongraphe}[4]{
{\centering #1} \\
\begin{picture}(0,0)
\put(0,0){\makebox(0,0)[b]{\rotatebox[origin=bl]{90}{#3}}}
\put(0,0){\makebox(0,0)[bl]{#4}}
\end{picture}{\centering #2}
}
\begin{document}

\doublespace
%
\begin{frontmatter}
\title{Enhanced sampling schemes for MCMC based blind Bernoulli-Gaussian deconvolution}

\author[ECN]{D. Ge\corauthref{cor}}
\corauth[cor]{Corresponding author. Tel.:+33 2 40 37 69 24; fax: +33 2 40 37 69 30}
\ead{Di.Ge@irccyn.ec-nantes.fr}
\author[ECN]{J. Idier}
\ead{Jerome.Idier@irccyn.ec-nantes.fr}
\author[ECN]{E. Le Carpentier}
\ead{Eric.Le-Carpentier@irccyn.ec-nantes.fr}

\address[ECN]{IRCCyN (CNRS UMR 6597), 1 rue de la Noë, BP 92101, 44321 Nantes Cedex 3, France}

\begin{abstract}
This paper proposes and compares two new sampling schemes for sparse deconvolution using a Bernoulli-Gaussian model. 
To tackle such a deconvolution problem in a blind and unsupervised context, the Markov Chain Monte Carlo (MCMC) framework is usually adopted, and the chosen sampling scheme is most often the Gibbs sampler. 
However, such a sampling scheme fails to explore the state space efficiently.
Our first alternative, the $K$-tuple Gibbs sampler, is simply a grouped Gibbs sampler. 
The second one, called partially marginalized sampler, is obtained by integrating the Gaussian amplitudes out of the target distribution.
While the mathematical validity of the first scheme is obvious as a particular instance of the Gibbs sampler, a more detailed analysis is provided to prove the validity of the second scheme.

For both methods, optimized implementations are proposed in terms of computation and storage cost. 
Finally, simulation results validate both schemes as more efficient in terms of convergence time compared with the plain Gibbs sampler. Benchmark sequence simulations show that the partially marginalized sampler takes fewer iterations to converge than the $K$-tuple Gibbs sampler. However, its computation load per iteration grows almost quadratically with respect to the data length, while it only grows linearly for the $K$-tuple Gibbs sampler.  
\end{abstract}

\end{frontmatter}


\section{Introduction}
\label{sec:intro}
The problem of the restoration of a sparse spike train \xb distorted by a linear system \hb and corrupted by noise $\epsilon$ such that $\zb = \xb \ast \hb + \epsilon$, arises in many fields such as seismic exploration~\cite{Cheng96,Kormylo83a} and astronomy~\cite{Bourguignon05}.

In this paper, we adopt a Bernoulli-Gaussian (BG) model for the spike train \xb, following~\cite{Kormylo83a} and many subsequent contributions such as~\cite{Cheng96,Champagnat96a}.
A BG signal is an independent, identically distributed (iid) process defined in two stages. Firstly, the sparse nature of the spikes is governed by the Bernoulli law:
\beq
\label{eqn:BG1}
P(\qb)=\lambda^L(1-\lambda)^{M-L}
\eeq
with the Bernoulli sequence $\qb=[q_1,\ldots,q_M]\T$ and $L=\sum_{m=1}^M q_m$, the number of non-zero realizations. Secondly, amplitudes $\xb=[x_1,\ldots,x_M]\T$ are assumed iid zero-mean Gaussian conditionally to \qb:
\beq
\label{eqn:BG2}
\xb\I\qb\ \sim \ \mathcal{N}(\zerob,\sigma_x^2\,\diag(\qb)),
\eeq
where $\diag(\qb)$ denotes a diagonal matrix whose diagonal is $\qb$. 

The MCMC approach~\cite{Liu01,Robert04} is a powerful numerical tool, appropriate to solve complex inference problems such as blind deconvolution. In the field of blind BG deconvolution, Cheng \etal\ pioneered the introduction of MCMC methods~\cite{Cheng96}. They proposed to rely on a \emph{plain} Gibbs sampler, \ie with a site-by-site updating scheme for the spike train, for which their algorithm constitutes a simple and canonical example.
However, simulation results indicate that it lacks reliability: from different initial conditions, significantly different estimations are obtained, even after a considerable number of iterations.

The recent contribution of~\cite{Labat06b} already identified a convergence issue linked to time-shift ambiguities, and proposed an efficient way to solve it. 
In addition, the scale ambiguities are treated by~\cite{Veit08} by proposing a scale re-sampling step to accelerate the convergence rate of the Markov chain.
In this study, we point out another source of inefficiency, unrelated to the above-mentioned ambiguities: instead of exploring the $2^M$ configurations of \qb at an acceptable speed, the Gibbs sampler tends to get stuck for many iterations around some particular configurations of \qb, often corresponding to local optimal configurations of \xb of the posterior distribution as illustrated by the example in Section~\ref{subsec:Classical Gibbs Sampling}. 
This conclusion meets Bourguignon and Carfantan's analysis~\cite{Bourguignon05}: the Markov chain equilibrates rapidly around a mode (\ie a local optimal configuration), but takes a long time to move from mode to mode.

In order to make up for this deficiency, our first proposition is to adopt a grouped Gibbs sampler~\cite{Liu01}, called a
\emph{K-tuple} Gibbs sampler, where blocks of $K$ adjacent BG variables $(q_i,\ldots,q_{i+K-1})$ and their associated amplitudes $(x_i,\ldots,x_{i+K-1})$ are jointly sampled.
A comparable idea first appeared in~\cite{Chi84}, and more recently in~\cite{Bourguignon05}, under the form of a deterministic iteration aimed at producing an increase of the posterior probability, whereas our goal is to sample the latter.

We then propose a second solution based on sampling the posterior distribution \emph{marginally} to the amplitudes \xb~\cite{Ge08}. A comparable idea is found in statistical signal segmentation~\cite{Dobigeon07} where some hyper-parameters are partially marginalized.
In fact, as Liu pointed out in~\cite{Liu01}, completely integrating out some components (the Gaussian amplitudes \xb in our case), leads to a more efficient sampling scheme called \emph{collapsed Gibbs sampling}. 
However, a plain collapsed Gibbs sampling on the marginal posterior distribution involves hardly tractable sampling steps. In particular, it is all but simple to sample \hb conditional on $(\zb,\qb)$ and marginally with respect to \xb. 
Our scheme solves this problem by combining a step that samples \qb marginally with respect to \xb and other sampling steps involving \xb.
Such a \emph{partially marginalized sampler} is fully valid from the mathematical point of view. In this paper, we show that it can be interpreted as a plain Gibbs sampler with a particular scanning order of the variables.  

Simulation tests on toy examples and on the so-called Mendel's sequence~\cite{Kormylo83a} 
confirm the efficiency of both methods in terms of computation time before convergence of the Markov chain. Further analysis shows the data length as a criterion to choose between the two proposed methods.  

This paper is organized as follows.  
In Section~\ref{sec:formulation}, after a brief formulation of the blind BG deconvolution problem, the Gibbs sampler of the joint posterior distribution~\cite{Cheng96} is presented, and an example illustrates its inefficiency as regards the sampling of \qb.
Sections ~\ref{sec:Generalized Gibbs on K-tuple variables} and~\ref{sec:Partially marginal Gibbs sampler} respectively introduce the generalized $K$-tuple Gibbs sampler and the partially marginalized sampler. In both cases, a toy example is used to evaluate the capability of the sampler to escape from local optimal configurations, and implementation issues are carefully dealt with.

Finally, simulation results are presented in Section~\ref{sec:Simulation results} to compare the efficiency of the sampling schemes according to Brooks and Gelman's convergence diagnostic~\cite{Brooks98}, and conclusions are drawn in Section~\ref{sec:Conclusion}.

\section{Problem formulation}
\label{sec:formulation}
\subsection{Statistical model}
\label{subsec:Statistical model}
The mathematical model of convolution reads
\beq
\label{eqn:direct}
z_n=\sum_{k=0}^Ph_k\,x_{n-k}+\epsilon_n
\eeq
for all $n \in \{1,\ldots, N\}$, where $\zb=\cro{z_1,\ldots,z_N}\T$ denotes the observed vector, $\xb=\cro{x_1,\ldots,x_M}\T$ is the unknown spike train, $\hb=\cro{h_0,\ldots,h_P}\T$ the impulse response (IR) of the system (assumed finite here) and $\epsilonb=\cro{\epsilon_1,\ldots,\epsilon_M}\T$ an noise vector, often assumed white stationary Gaussian. 
The deconvolution problem is said \emph{blind} when $\hb$ is unknown, which is the studied case here.
Akin to~\cite{Cheng96}, the following assumptions are made:
\bit
\item $\epsilonb\sim\mathcal{N}(0,\sigma_{\epsilon}^2\Iv)$ is independent of \xb and \hb;
\item $\xb$ is a BG process defined by Eqs.\eqref{eqn:BG1}-\eqref{eqn:BG2} with $\sigma_{\xb}^2=1$;
\item $\hb \sim \mathcal{N}(0,\sigma_h^2\Iv_{P+1})$.
\eit
Let us remark that by setting $\sigma_{\xb}^2$ arbitrarily to one while imposing the Gaussian prior law \hb removes the scale ambiguity inherent to the blind deconvolution problem~\cite{Labat06b}. 
By adopting the ``zero boundary'' condition so that $M=N-P$, \zb can be rewritten in the form of a matrix multiplication $\zb = \Hv\xb + \epsilonb$, where \Hv denote the $N\times M$ Toeplitz matrix of convolution. 
It is also useful to introduce the Toeplitz matrix \Xv of size $N\times (P+1)$ such that $\Hv\xb = \Xv\hb$. 
According to the Monte Carlo principle, a posterior mean estimator of $\Theta = \{\qb,\xb,\hb,\lambda,\sigma^2_{\epsilon},\sigma^2_{\hb}\}$ given \zb can be approximated by:
\begin{equation}\label{eqn:MonteCarlo}
\widehat{\Theta} = \frac1{I-J}\sum_{k=J+1}^I\Theta^{(k)},
\end{equation}
where the sum extends over the last $I-J$ samples. In the MCMC framework, the samples are generated iteratively, so that asymptotically  $\Theta^{(k)}$ follows the joint posterior distribution~\cite{Robert04}:
\bal
P(\Theta\,|\,\zb) \propto~
&g(\zb-\Hv\xb;\sigma_{\epsilon}^2\Iv_N)\, 
g(\xb;\sigma_{\xb}^2\text{diag}\{\qb\})\,
g(\hb;\sigma_h^2\Iv_{P+1})\,P(\qb;\lambda)\,P(\sigma^2_{\hb})\,P(\lambda)\,P(\sigma^2_{\epsilon})
\label{eqn:JointPostdistribution}
\eal
where $g(\cdot;\Rv)$ denotes the centered Gaussian density of covariance \Rv. Conjugate prior laws are adopted for the last three terms $P(\sigma^2_{\hb})$,$P(\lambda)$ and $P(\sigma^2_{\epsilon})$:  
\balx
\sigma_{\hb}^2 \sim \mathit{IG}(1,1),\quad \lambda \sim \mathit{Be}(1,1), \quad \sigma_{\epsilon}^2 \sim \mathit{IG}(1,1),
\ealx
where $\mathit{IG}$ and $\mathit{Be}$ respectively represent the inverse Gamma and Beta distributions. The flat shapes of the two distributions $\mathit{IG}(1,1)$ and $\mathit{Be}(1,1)$ (a uniform distribution on $[0,1]$) can be interpreted as conveying no specific prior information of the parameters. 

\subsection{Classical Gibbs Sampling}
\label{subsec:Classical Gibbs Sampling}
A Gibbs sampler circumvents the difficulties of directly inferencing from the joint posterior distribution as in Eq.~\eqref{eqn:JointPostdistribution} and instead draws each parameter from its conditional law while all others are fixed. The principle of Cheng \etal's Gibbs sampler is given in Table~\ref{Cheng0} and a pseudocode can be found in Table~\ref{Cheng}.
\begin{table}[b]
\caption{Cheng \etal's Gibbs sampler.}
\label{Cheng0}
\centering
\medskip
\bfmi[.45\linewidth]
\btabu{ll}
\textcircled{\tiny1} &Let $\yb =(\qb,\xb)$. For each $i=1\ldots,M,$\\
& draw $y_i^{(k+1)}\I\yb_{1:i-1}^{(k+1)},\yb_{i+1:M}^{(k)},\hb^{(k)},\sigma_{\epsilon}^{(k)},\lambda^{(k)}, \zb$\\ 
\textcircled{\tiny2} & draw $\hb^{(k+1)} \I\xb^{(k+1)},\sigma_{\hb}^{(k)},\sigma_{\epsilon}^{(k)}, \zb$ \\
\textcircled{\tiny3} & draw $\sigma_{\epsilon}^{(k+1)}\I\xb^{(k+1)},\hb^{(k+1)}, \zb$\\ 
\textcircled{\tiny4} & draw $\lambda^{(k+1)}\I\qb^{(k+1)}$\\
\textcircled{\tiny5} & draw $\sigma_{\hb}^{(k+1)}\I\hb^{(k+1)}$ 
\etabu
\efmi
\end{table}
\begin{table}[b]
\caption{A pseudocode for Cheng \etal's Gibbs sampler (see~\cite{Cheng96} for implementation details).
}
\label{Cheng}
\centering
\medskip
\bfmi[.8\linewidth]
\begin{algorithmic}
\STATE \COMMENT{Initialization}
\STATE $\hb\gets\zerob;~\hb(\texttt{round}(P/2))\gets1$ 
\STATE $\qb\gets\zerob;~\xb\gets\zerob$ 
\STATE $\sigma_{\xb}^2\gets1$
\STATE Sample $\sigma_{\epsilon}^2 \sim \mathit{IG}(1,1);~\sigma_{\hb}^2 \sim \mathit{IG}(1,1)$
\REPEAT 
\STATE \COMMENT{--------------------- Step~1: Sample $y_i = (q_i,x_i)$ ---------------------}
\STATE $\sigma_1^2\ \gets \sigma_{\epsilon}^2\sigma_{\xb}^2 / (\sigma_{\epsilon}^2 + \sigma_{\xb}^2\left\|\hb\right\|^2)$ 
\STATE $\eb\gets\zb-\Hv\xb$
\FOR{$i=1$ to $M$}
\STATE $\eb_i\gets\eb+\hb_i x_i$ \COMMENT{$\hb_i$ is the $i$-th column of \Hv}
\STATE $\mu_i \gets (\sigma_1^2 / \sigma_{\epsilon}^2) \hb_i\T\eb_i$
\STATE ${\nu_i} \gets \lambda (\sigma_1 / \sigma_{\xb})\exp\left(\mu_i^2 / (2\sigma_1^2)\right)$
\STATE $\lambda_i \gets {\nu_i}/({\nu_i}+1-\lambda)$
\STATE Sample $q_i \sim \textit{Bi}(\lambda_i)$ \COMMENT{Step 1a}
\STATE Sample $x_i \sim \Nc(\mu_i,\sigma_1^2)$ if $q_i=1$, $x_i=0$ otherwise \COMMENT{Step1b} 
\ENDFOR
\STATE \COMMENT{--------------------- Step~2: Sample $\hb$ ---------------------}
\STATE $\Sv \gets \sigma_{\epsilon}^{-2}\Xv\T\Xv+\sigma_{\hb}^{-2}\texttt{eye}(P+1)$ \COMMENT{Inverse of covariance matrix $\Rv$} 
\STATE $\Uv \gets \texttt{chol}(\Sv)$ \COMMENT{Cholesky factor of $\Sv=\Rv\M$, \ie $\Rv=\Uv\M\Uv\MT$}
\STATE $\hb \gets \Uv\backslash \left(\sigma_{\epsilon}^{-2} \Uv\T\backslash \Xv\T\zb + \texttt{randn}(P+1,1)\right)$ \quad \COMMENT{$\hb \sim \Nc(\mb,\Rv)$ with $\mb=\sigma_{\epsilon}^{-2}\Rv\Xv\T\zb$}
\STATE \COMMENT{--------------------- Step~3: Sample $\sigma_{\epsilon}^2$ ---------------------}
\STATE Sample $\sigma_{\epsilon}^2 \sim \text{IG}(N/2+1,\left\|\zb-\Hv\xb\right\|^2/2+1)$
\STATE \COMMENT{--------------------- Step~4: Sample $\lambda$ ------------------------}
\STATE Sample $\lambda \sim \text{Be}(1+L,1+M-L)$, with $L=\sum_m q_m$
\STATE \COMMENT{--------------------- Step~5: Sample $\sigma_{\hb}^2$ ---------------------}
\STATE Sample $\sigma_{\hb}^2 \sim \text{IG}(P/2+1,\left\|\hb\right\|^2/2+1)$
\UNTIL{Convergence}
\end{algorithmic}
\efmi\\
\end{table}
The five simulation steps are iterated until convergence toward the posterior distribution in Eq.~\eqref{eqn:JointPostdistribution}, and $\widehat{\Theta}$ is finally built according to (\ref{eqn:MonteCarlo}).
Let us remark than in Table~\ref{Cheng}, Step 1 is performed in two substeps: for each $i$, $q_i$ is first sampled conditional on $\left( \Theta \setminus (q_i, x_i), \zb \right)$, and then $x_i$ conditional on $\left( \Theta \setminus x_i, \zb\right)$. Taken jointly, these two substeps perform the sampling of $y_i=(q_i, x_i)$ conditional on $\left( \Theta \setminus y_i, \zb \right)$.

Despite its simplicity, this scheme presents several drawbacks.
In~\cite{Labat06b}, it is pointed out that time-shift ambiguities may lead to unreliable estimates, depending on the initialization of \hb.
To make up for this deficiency, an additional Metropolis-Hasting (MH) step is proposed between Steps~1 and~2 in order to allow moves to a time-shifted version of $(\yb,\hb)$. 
Another inefficiency in the Cheng \etal's Gibbs sampler has been discussed in~\cite{Veit08} with regard to the scale factor between \hb and \xb. A scale re-sampling step is proposed and tested to yield faster convergence. 
In Table~\ref{algo:TimeShift}, the modified version of Step~2 incorporates the time-shift operation and the scale re-sampling with an efficient implementation in which $f_{\text{GIG}}$ represents the \emph{Generalized Inverse Gaussian} law (see~\cite{Labat06b} and~\cite{Veit08} for more details). 
In what follows, the resulting global scheme is referred to as a hybrid sampler, since there is a Metropolis Hastings step within the Gibbs scheme and no more a pure Gibbs sampler. 

\begin{table}[h]
\caption{Modified version of Step~2 to incorporate the time-shift compensation according to~\cite{Labat06b} and the scale re-sampling according to~\cite{Veit08} .
}
\label{algo:TimeShift}
\centering
\medskip
\bfmi[.8\linewidth]
\begin{algorithmic}
\STATE \COMMENT{ --------------------- Time-shift compensation (Labat \etal~\cite{Labat06b})------------}
\STATE Propose $\yb'=(\qb',\xb')$ from $\yb=(\qb,\xb)$ with probability\\[1mm]
\centerline{$\pi(\yb'\I \yb)=\left\{
\barr{ll}
1-2\eta &\text{if}\; \yb' = \yb,\\
\eta &\text{if}\; \yb' = (\texttt{circshift}(\qb,1),\texttt{circshift}(\xb,1)),\\
\eta &\text{if}\; \yb' = (\texttt{circshift}(\qb,-1),\texttt{circshift}(\xb,-1)), \\
\earr
\right.
$}\smallskip
where $\eta\in(0,1/2)$, \eg $\eta=1/4$.
\STATE $\Sv \gets \sigma_{\epsilon}^{-2}\Xv\T\Xv+\sigma_{\hb}^{-2}\texttt{eye}(P+1)$ \COMMENT{Inverse of covariance matrix $\Rv$} 
\STATE $\Sv' \gets \sigma_{\epsilon}^{-2}(\Xv')\T\Xv'+\sigma_{\hb}^{-2}\texttt{eye}(P+1)$ \COMMENT{Inverse of covariance matrix $\Rv$} 
\STATE $\Uv \gets \texttt{chol}(\Sv)$
\STATE $\Uv' \gets \texttt{chol}(\Sv')$
\STATE $\rho \gets \sigma_{\epsilon}^{-4}\norm{(\Uv')\T\backslash (\Xv')\T\zb}^2 -\sigma_{\epsilon}^{-4}\norm{\Uv\T\backslash \Xv\T\zb}^2 + 2\texttt{log} |\Uv|/|\Uv'|$
\STATE \COMMENT{$\rho = (\mb')\T(\Rv')\M \mb' - \mb\T\Rv\M \mb + \log |\Rv'|/|\Rv|$}
\IF[with probability $\min\{1,\exp(\rho/2)\}$]{$2\texttt{log(rand)}<\rho$}
\STATE $\yb \gets \yb'$
\STATE $\hb \gets \Uv'\backslash \left(\sigma_{\epsilon}^{-2} (\Uv')\T\backslash (\Xv')\T\zb + \texttt{randn}(P+1,1)\right)$
\ELSE
\STATE $\hb \gets \Uv\backslash \left(\sigma_{\epsilon}^{-2} \Uv\T\backslash \Xv\T\zb + \texttt{randn}(P+1,1)\right)$
\ENDIF
\STATE \COMMENT{ --------------------- Scale re-sampling (Veit \etal~\cite{Veit08})------------}
\STATE $\lambda \leftarrow (\sum q_i -P-1)/2$
\STATE $\alpha \leftarrow |\xb|^2/\sigma_{\xb}^2, \quad \beta \leftarrow |\hb|^2/\sigma_{\hb}^2$
\STATE $s^2 \sim f_{\text{GIG}}(\lambda, \alpha, \beta)$ \quad \COMMENT{simulation by acceptation-rejection}
\STATE $\xb \leftarrow \xb\cdot s, \quad \hb \leftarrow \hb / s$
\end{algorithmic}
\efmi
\end{table}

\begin{figure}[htb]
  \centering
\btabu{@{}c@{~~}c@{}}
  \figc[width=.47\linewidth]{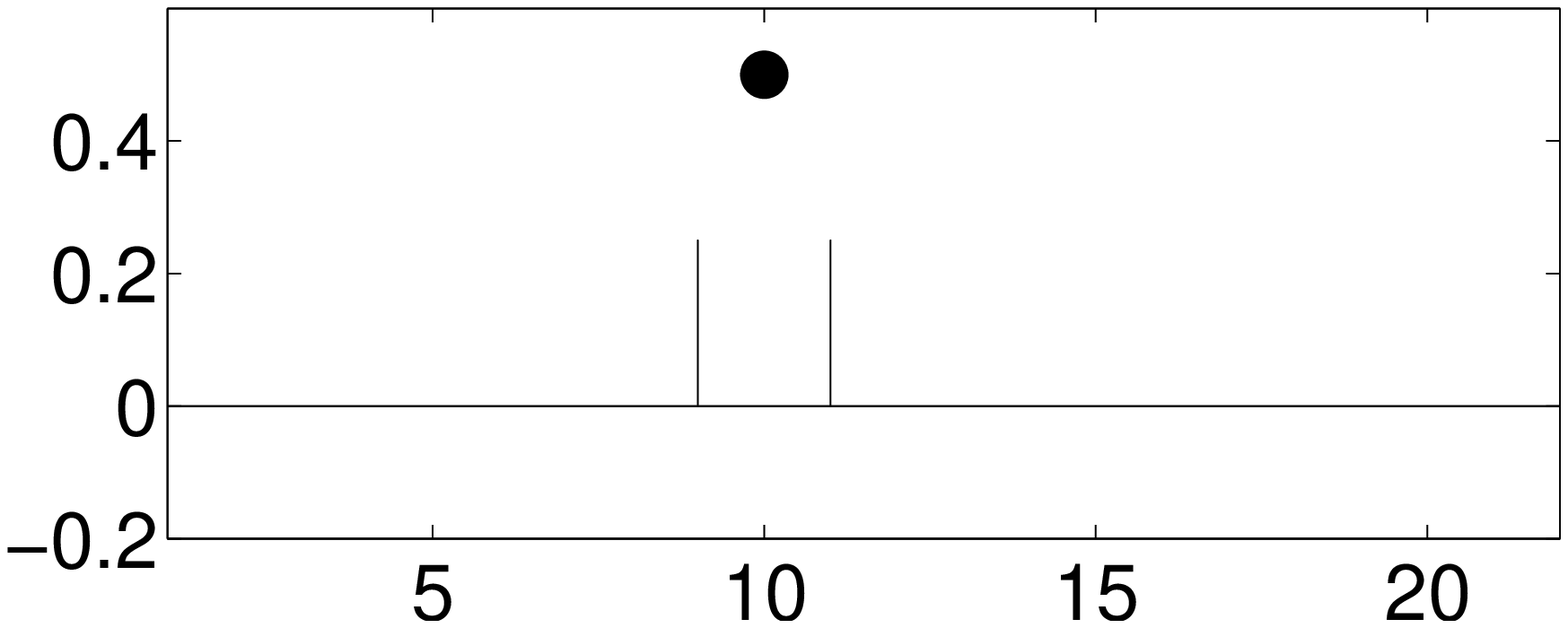} & \figc[width=.47\linewidth]{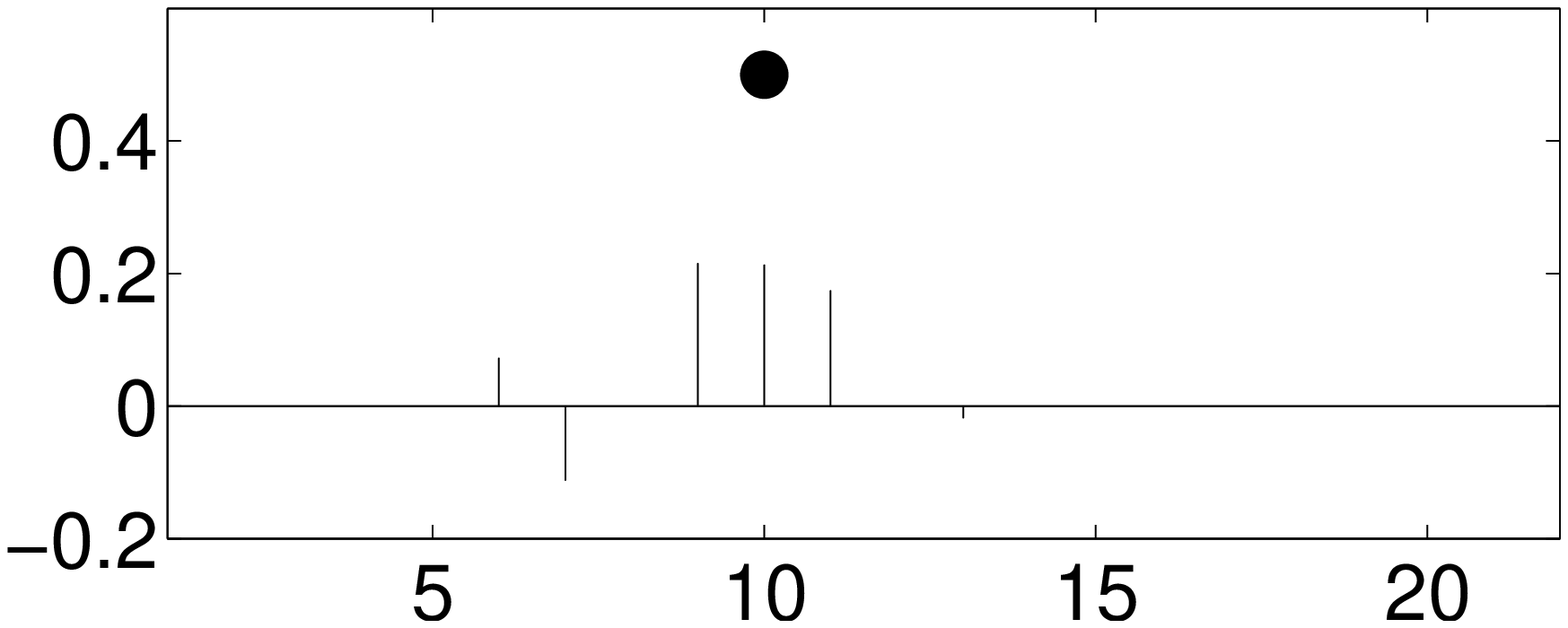} \\
{\footnotesize  (a) Initial configuration.} & {\footnotesize (b) 1000th iteration}\\[2mm]
  \figc[width=.47\linewidth]{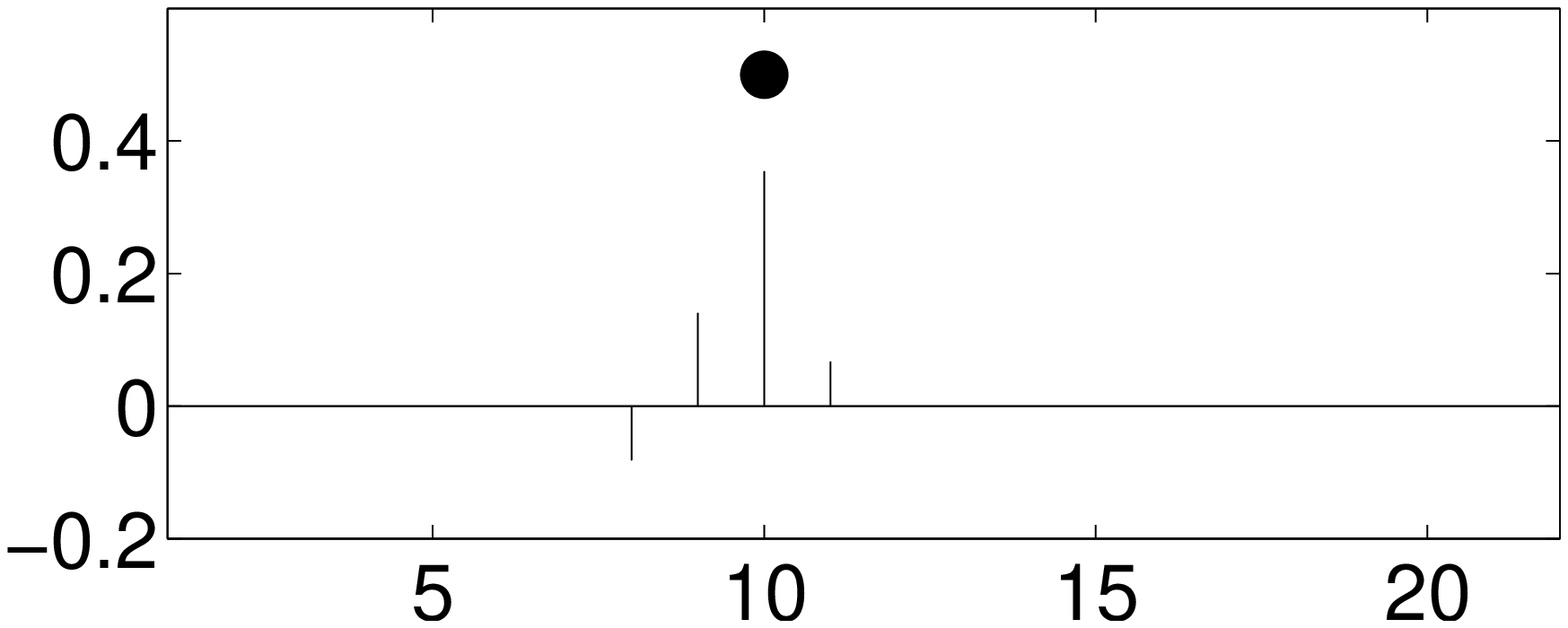} & \figc[width=.47\linewidth]{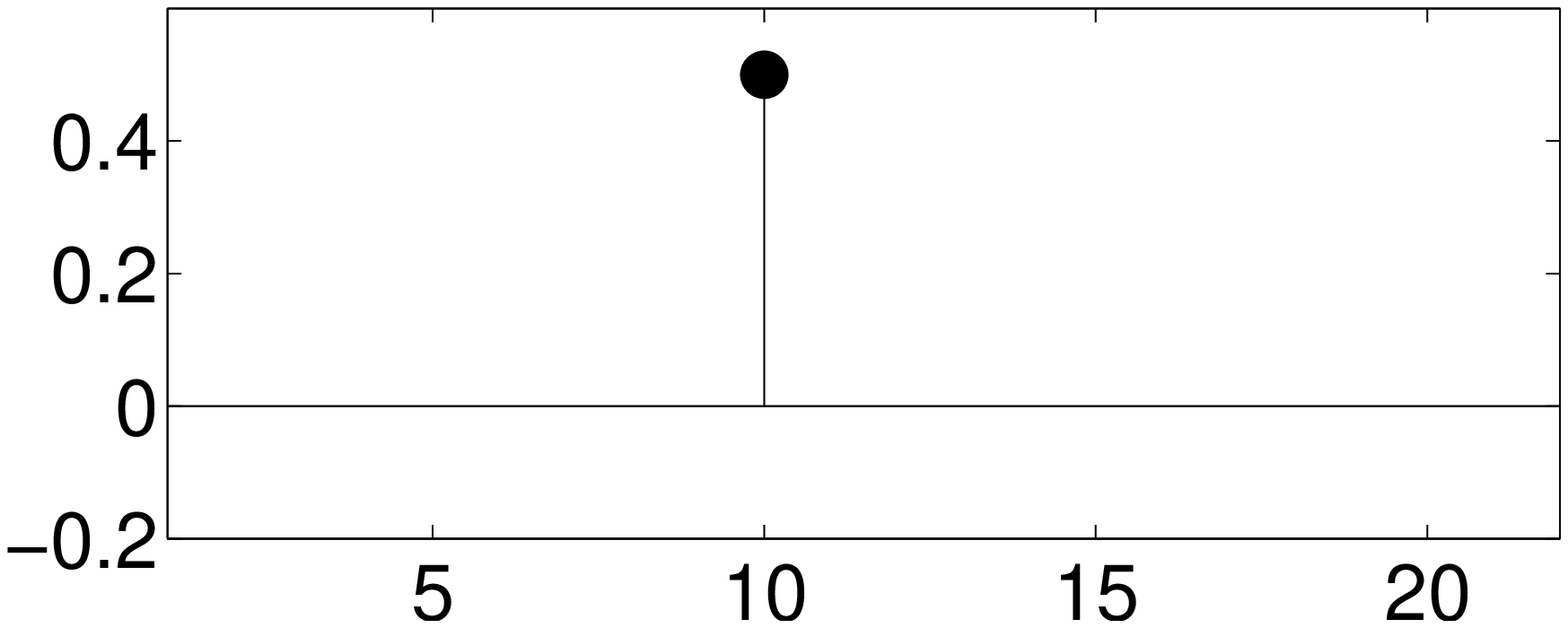} \\
{\footnotesize  (c)  2000th iteration} & {\footnotesize (d) 2415th iteration}\\
\etabu
\caption{Sampled BG sequences obtained by the hybrid sampler on a simple example. From a local optimal configuration chosen as initial state, 
the Markov chain spends several hundreds of iterations before visiting the solution, \ie a unique spike in position 10 (marked as a bullet). 
}
\label{fig:Configs}
\end{figure}

\begin{figure}[htb]
  \centering
\setlength{\tabcolsep}{0pt}
\btabu{cc}
  \figc[height=6cm]{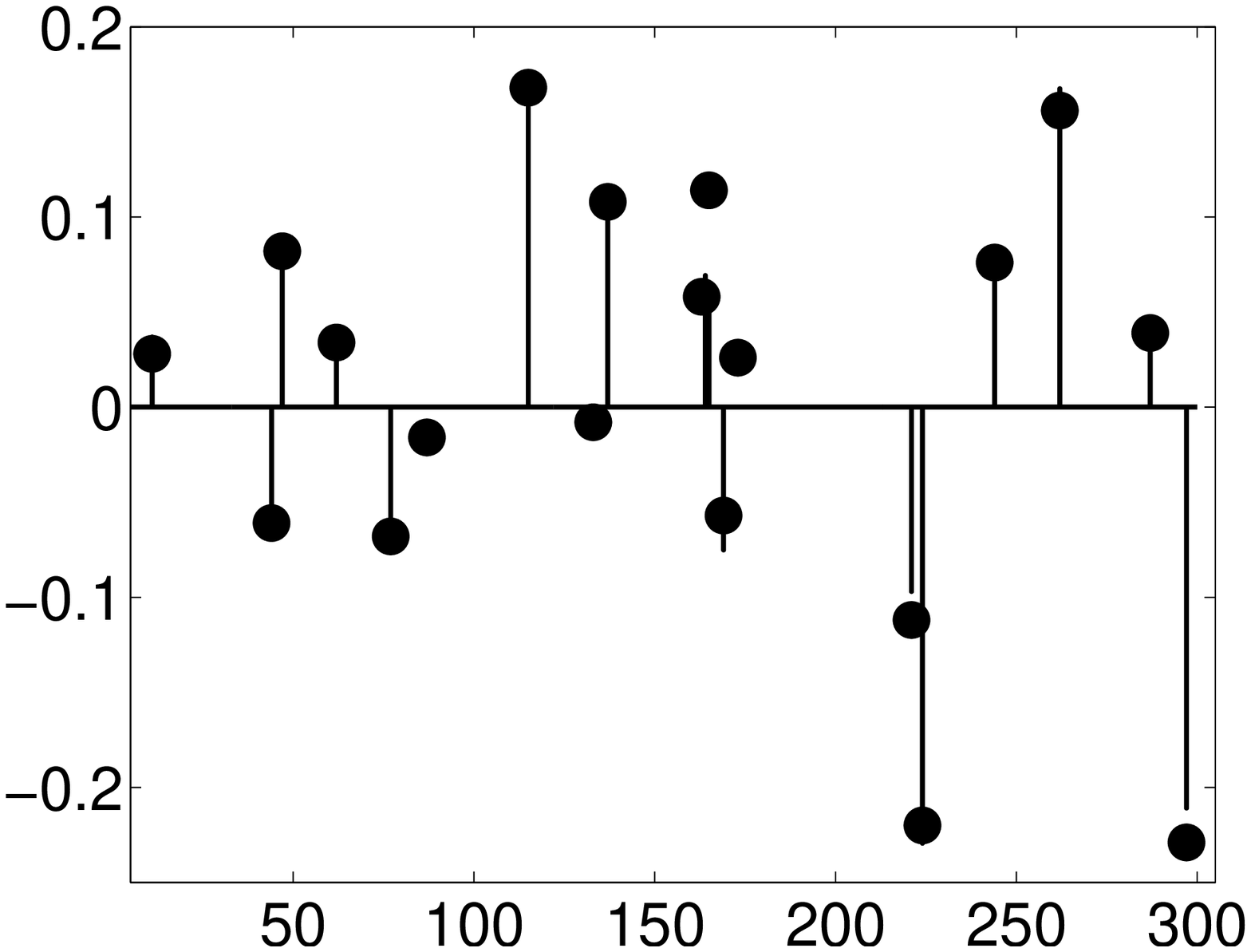}
&  \figc[height=6cm]{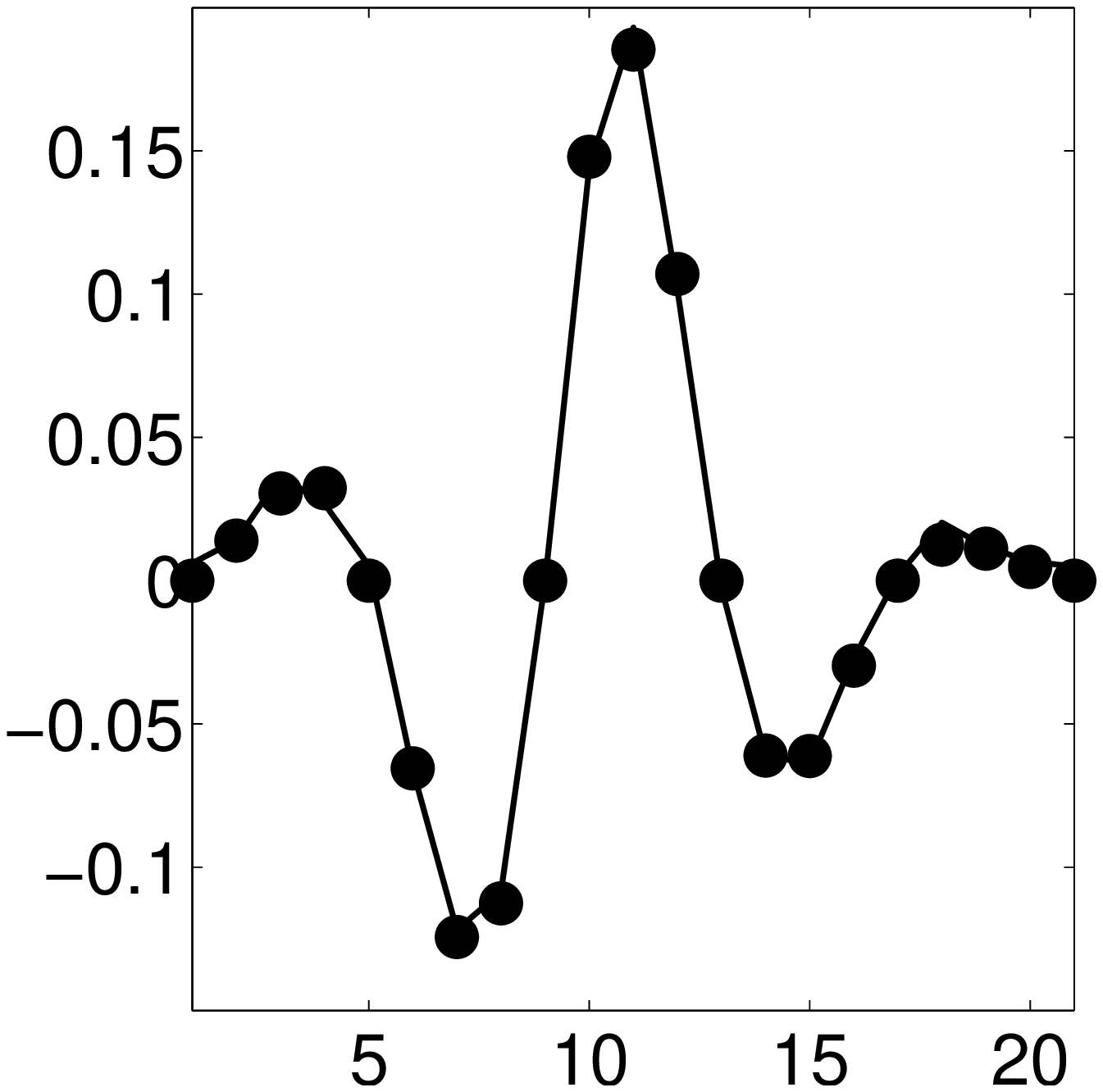}
\\
  \figc[height=6cm]{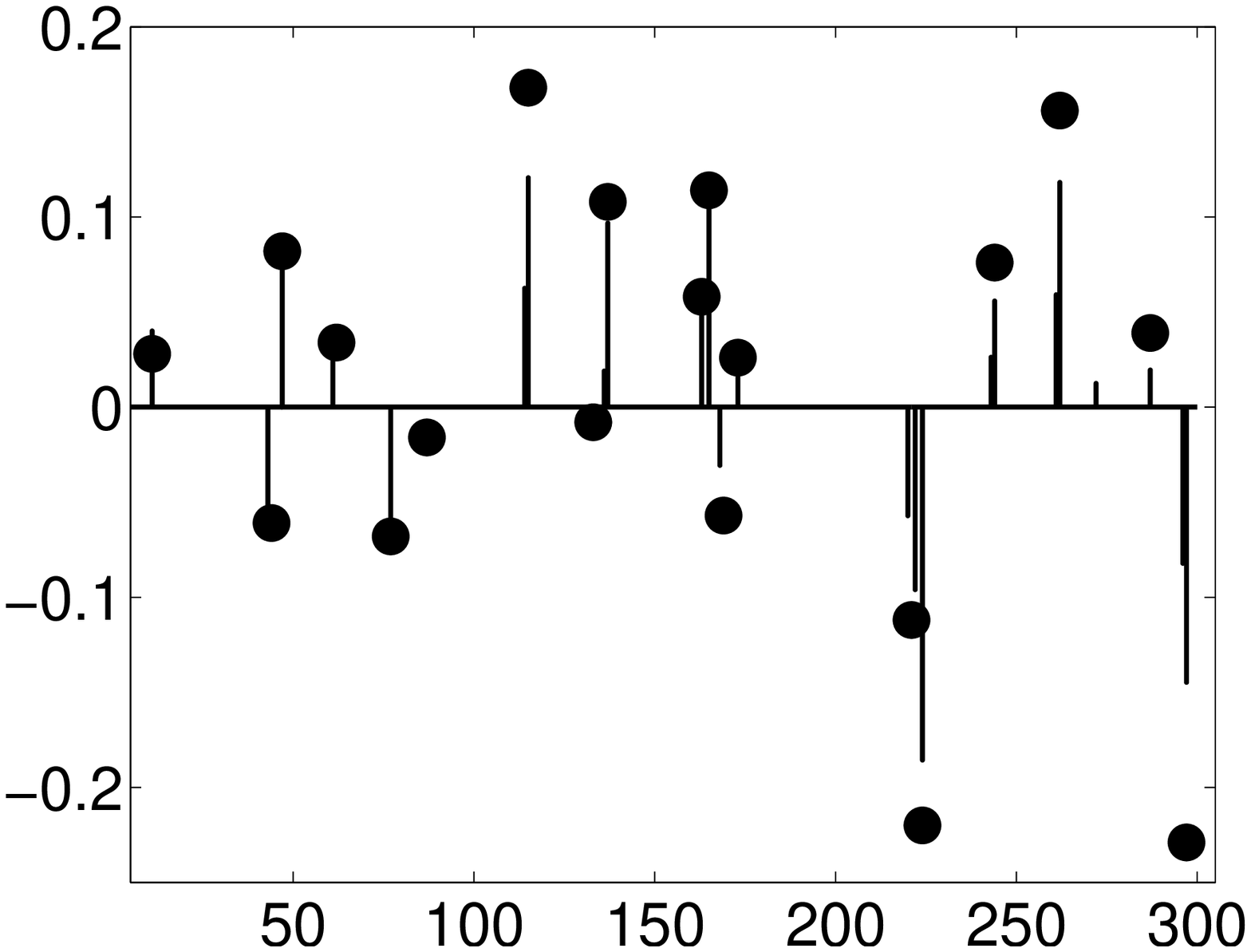}
&  \figc[height=6cm]{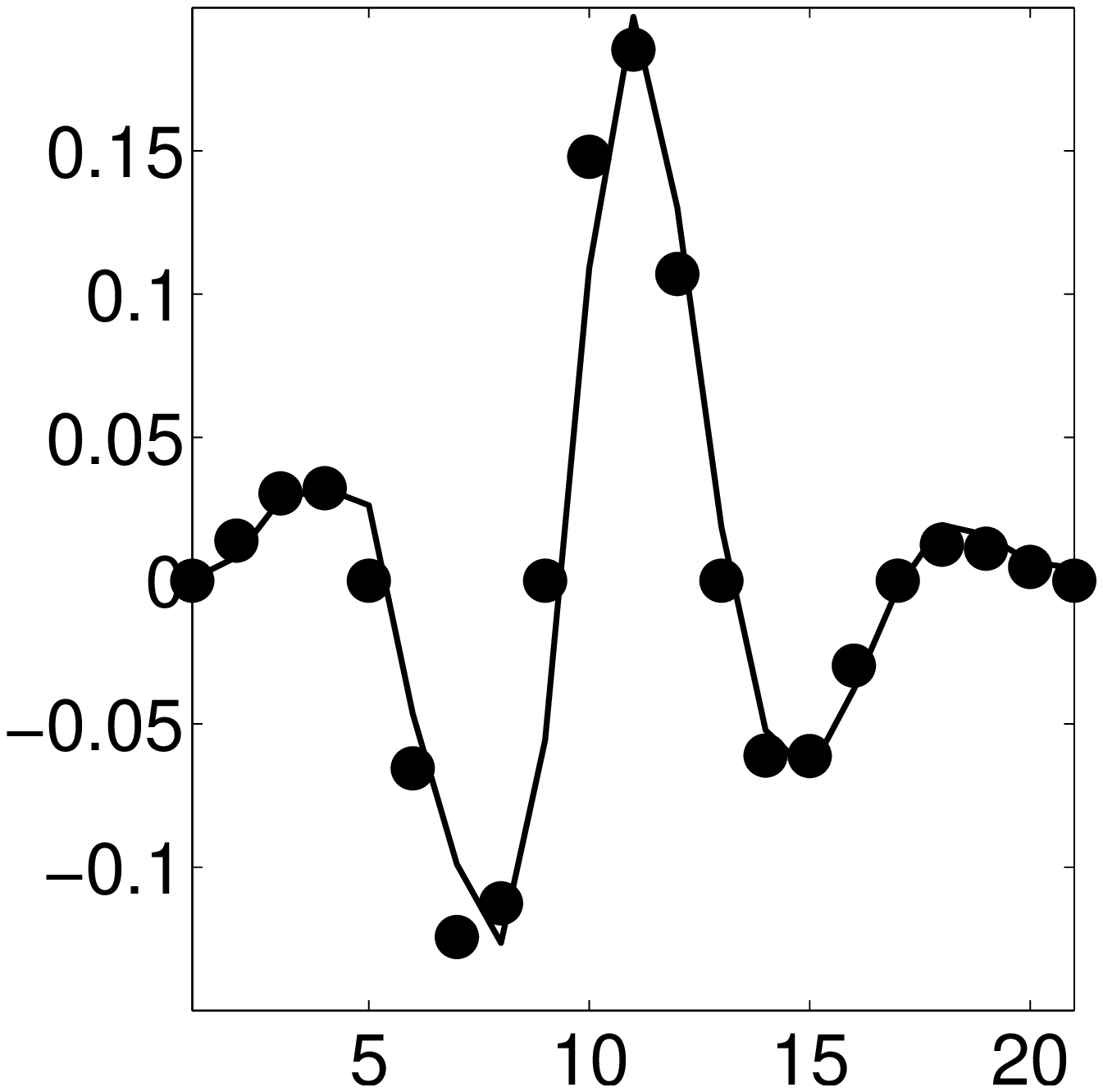}
\\
  \figc[height=6cm]{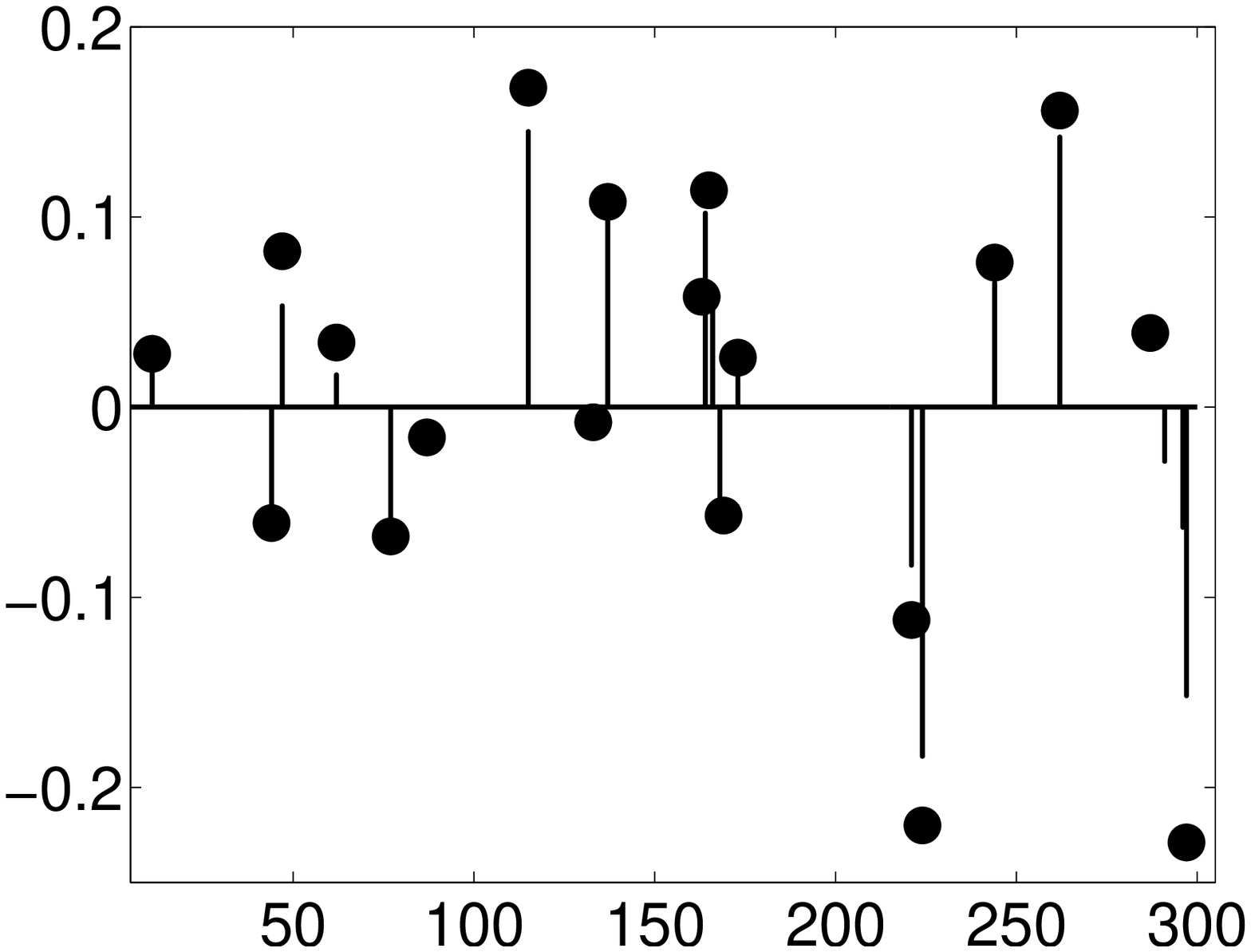}
&
  \figc[height=6cm]{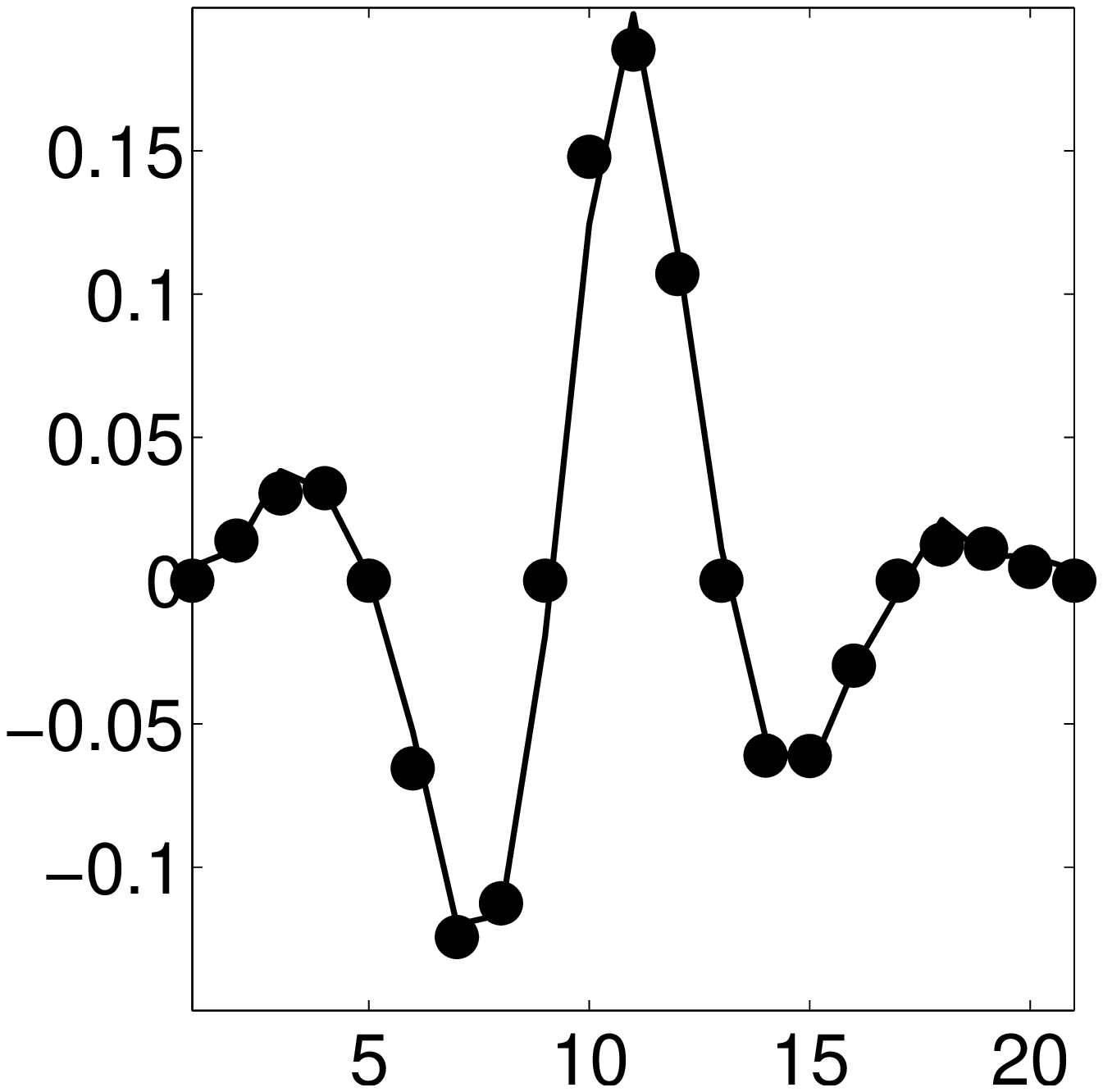}
\\
  {\footnotesize (a) estimation of \xb}
&  {\footnotesize (b) estimation of \hb}
\etabu
\caption{Three different solutions of the hybrid sampler~\cite{Labat06b}, obtained by changing the random generator seed only. The actual values are marked as bullets.}
\label{fig:MCMCTimeshifted}
\end{figure}

Despite the improvement brought by the timeshift and scale ressampling operations in the hybrid sampler, we have identified another source of inefficiency, completely independent from the time-shift and/or scale ambiguity problem. Figure~\ref{fig:Configs} is an illustration based on simulated data generated from a single spike, convolved with an IR defined by 
$$
\hb(i)=\cos\pth{(i-10)\frac{\pi}{4}}\exp\pth{-|0.225i-2|^{1.5}},\ i=0,\ldots,20,
$$
which is depicted by bullets on Figure~\ref{fig:MCMCTimeshifted}(b). A local optimal configuration of \xb (two neighboring spikes instead of a single one) is chosen as the initial state in Figure~\ref{fig:Configs}(a).
The hybrid sampler generates a Markov chain that typically takes several thousands of iterations (depending on the chosen random generator seed) to visit the optimal configuration for the first time.
%
Figure~\ref{fig:MCMCTimeshifted} further illustrates that the hybrid sampler tends to produce unreliable estimated values. Here, Mendel's well-known BG sequence is adopted as a test signal~\cite{Kormylo83a}. The same IR is used 
and the data are corrupted by Gaussian noise with $\sigma_{\epsilon}^2=4\times10^{-6}$, corresponding to a signal-to-noise ratio (SNR) of $12.80$\,dB. 
The three estimation results in Figure~\ref{fig:MCMCTimeshifted} are obtained from the same simulated data \zb, and the same initial states, the only difference being the seed value of the random generator. 
For each Markov chain, $1000$ samples are produced and the last 250 are averaged to compute the estimation.
Substantial variations exist from one estimated result to another, especially in the number and the positions of spikes in Figure~\ref{fig:MCMCTimeshifted}(a). Actually, it can be checked that each sequence $\stdacc{\qb^{(k)}}$ tends to remain constant for many iterations, instead of exploring the state space efficiently. It should be stressed that these results are typical, \ie they have not been selected on purpose.

This phenomenon is basically due to Step~1 of the sampler. 
More precisely, the corresponding sampled chains hardly escape from local optima of the posterior distribution, because the consecutive samples are highly dependent. 



\section{Generalized Gibbs on K-tuple variables}

\label{sec:Generalized Gibbs on K-tuple variables}
This inefficiency in MCMC sampling schemes has been already noticed by Bourguignon and Carfantan in~\cite{Bourguignon05}. They propose a solution involving shifts of detected spikes to adjacent positions.
However, such a solution is a deterministic procedure aimed at increasing the posterior probability during the burn-in period of the Markov chain.
Therefore, it does not leave the posterior distribution invariant. In contrast, our goal here is to propose a valid sampling step that produces spike shifts.

Let us introduce the set $\Omega=\{0,\ldots,K-1\}$, and the notation $i+\omega= \{i+m \I m \in \omega\}$ for all $\omega\subset\Omega$. For any vector \vb, let us denote $\vb_\omega$ the subvector formed by the entries $v_m$, $m\in\omega$, and $\vb_{-\omega}$ the subvector formed by the remaining entries. Similarly, for any matrix \Vb, we introduce the submatrices $\Vb_\omega$ and $\Vb_{-\omega}$ formed by columns of \Vb, indexed by $\omega$ in the first case, and gathering the remaining columns in the second case.

The principle of the grouped Gibbs sampler to K-tuple variables consists in updating vector $\yb_{i+\Omega}$ for each $i=1,\ldots,M-K+1$ instead of the scalar updates of Step 1 in Table~\ref{Cheng0}.
%
In particular, let us remark that such a joint update strategy allows a permutation of the Bernoulli vector $\qb_{i+\Omega}$ (from a configuration $(0,1)$ to $(1,0)$ in the case of $K=2$) within one iteration.

The resulting sampling scheme is given in Table~\ref{Generalized Gibbs}. 
The notation of $\qb_{i+\Omega}^{(k+1/2)}$ (and that of $\xb_{i+\Omega}^{(k+1/2)}$) denotes the fact that $\yb^{(k+1)}$ is latter re-sampled in step 2(a). 
Compared to the hybrid version of the original Gibbs sampler of Table~\ref{Cheng0}, only Step~1 has been modified. It is clear that in the case of $K=1$, we are driven back to the hybrid sampler.
\begin{table}[htb]
\caption{Grouped Gibbs sampler updating $K$-tuple variables}
\label{Generalized Gibbs}
\centering
\medskip
\bfmi[.5\linewidth]
\btabu{ll}
\textcircled{\tiny1}&For each  $i=1\ldots,M-K+1$,  \\
&(a)~draw $\qb_{i+\Omega}^{(k+1/2)}\I \yb_{-(i+\Omega)}^{(k)},\hb^{(k)},\sigma_{\epsilon}^{(k)},\lambda^{(k)}, \zb$\\
&(b)~draw $\xb_{i+\Omega}^{(k+1/2)}\I \yb_{-(i+\Omega)}^{(k)},\qb_{i+\Omega}^{(k+1/2)},\hb^{(k)},\sigma_{\epsilon}^{(k)}, \zb$\\
\textcircled{\tiny2}&(a) draw $\yb^{(k+1)}$ according to Table~\ref{algo:TimeShift}\\
&(b) draw $\hb^{(k+1)}\I\yb^{(k+1)},\sigma_{\hb}^{(k)},\sigma_{\epsilon}^{(k)}, \zb$\\
\textcircled{\tiny3}& draw $\sigma_{\epsilon}^{(k+1)}\I\xb^{(k+1)},\hb^{(k+1)}, \zb$\\
\textcircled{\tiny4}& draw $\lambda^{(k+1)}\I\qb^{(k+1)}$\\
\textcircled{\tiny5}& draw $\sigma_{\hb}^{(k+1)}\I\hb^{(k+1)}$
\etabu
\efmi
\end{table}

It should be noted that the complexity of the $K$-tuple sampler increases rapidly (indeed, exponentially) with $K$, while larger values of $K$ are expected to allow a more efficient exploration of the state space.
The latter fact is illustrated by Table~\ref{tab:K-tupleIter}, which reports the number of iterations needed to visit the exact configuration of \qb for the first time, for different values of $K$ and of the random generator seed. For $K=1$, the required number of iterations varies in large amounts. Typically, several thousands of iterations are required. In contrast, the exact configuration of \qb is visited after a few tens of iterations for $K=2$. For $K=3$ or $K=4$, the first time visit happens immediatly in most cases.
\btabl[htb]
\caption{Number of iterations for $K$-tuple samplers to visit the true configuration for the first time in the example of Figure~\ref{fig:Configs}\label{tab:K-tupleIter}}
\centering
\medskip
\setlength{\tabcolsep}{3pt}
\btabu{|c|c|c|c|c|}
\hline
   \text{seed value}& $K=1$ & $K=2$ & $K=3$ & $K=4$ \\
\hline
$19$ &	2415 & 21  & 4 & 1 \\
\hline
$29$ &  1558 & 55 & 1 & 1   \\
\hline
$39$&  5137 & 12 & 1 & 1 \\
\hline
$49$&  7132 & 37 & 3 & 2 \\
\hline
\etabu
\etabl

In practice, the question remains to know what is the best trade-off between $K=1$ or larger values of $K$. In the first case, the convergence is slow in terms of iteration number, but the computing time for each iteration is relatively low. By increasing values of $K$, the computing time per iteration increases rapidly, while the iteration number decreases. Finding the best trade-off on a purely theoretical and general basis is probably hopeless. More modestly, we examine this question in Section~\ref{sec:Simulation results} on the sole basis of Mendel's example of Figure~\ref{fig:MCMCTimeshifted}. However, a prerequisite is to rely on an optimized implementation of the K-tuple algorithm. This is the goal of the next subsection.



\subsection{Numerical implementation}

Steps~2 to 5 being identical to those of the hybrid sampler, only the algorithmic implementation of the Step~1 is detailed in Table~\ref{algo:Ktuple}.
In the K-tuple sampler, Step~1 draws $(\qb_{i+\Omega},\xb_{i+\Omega})$ according to their joint conditional law:
\bal
P(\qb_{i+\Omega},\xb_{i+\Omega}\I \text{rest}) &\propto P(\zb\I \xb_{i+\Omega}, \xb_{-(i+\Omega)},\hb, \sigma_{\epsilon}^2) \, 
P(\xb_{i+\Omega}\I \qb_{i+\Omega},\sigma_x^2) \, P(\qb_{i+\Omega}\I \lambda)\notag\\
&\propto\exp\acc{-\frac1{2\sigma_{\epsilon}^2}\left\|\eb_{K,i}- \Hv_{i+\Omega}\xb_{i+\Omega}\right\|^2}
\, g(\xb_{i+\Omega};\sigma_x^2\diag(\qb_{i+\Omega})) \, 
\pth{\frac{\lambda}{1-\lambda}}^{\unb\T\qb_{i+\Omega}}
\label{jointlaw}
\eal
where 
$\eb_{K,i}= \zb- \Hv_{-(i+\Omega)}\xb_{-(i+\Omega)}$,
so that $\eb_{K,i}$ is a function of $(\xb_{-(i+\Omega)},\hb,\zb)$ but not of $\xb_{i+\Omega}$.

Let $\omega\subset\Omega$ be such that $\qb_{i+\omega}$ gathers the nonzero entries of $\qb_{i+\Omega}$, so that $\unb\T\qb_{i+\Omega}=\#\omega$.
From~\eqref{jointlaw}, it is straightforward to deduce that Step~1(b) amounts to draw a Gaussian law defined by $\xb_{i+\Omega\setminus\omega}=\zerob$ and, if $\omega\neq\emptyset,$ $\xb_{i+\omega}\sim\Nc(\mb_{i+\omega},\Rv_{i+\omega})$,
where
\bal
\label{R}
\Rv_{i+\omega}\M&={\sigma_{\epsilon}^{-2}}\Hv_{i+\omega}\T \Hv_{i+\omega} + {\sigma_{\xb}^{-2}} \Iv_{\#\omega},
\\
\label{m}
\mb_{i+\omega}&={\sigma_{\epsilon}^{-2}}\Rv_{i+\omega}\Hv_{i+\omega}\T\eb_{K,i}.
\eal

By integration of~\eqref{jointlaw} with respect to $\xb_{i+\Omega}$, it is then possible to express the marginal conditional law of $\qb_{i+\Omega}$, along which the latter vector must be sampled according to Step~1(a):
$$
P(\qb_{i+\omega}=\unb,\qb_{i+\Omega\setminus\omega}=\zerob) = \frac{p_{i+\omega}}{\sum_{\omega'}p_{i+\omega'}}
$$
where
\beq
\label{p}
p_{i+\omega}=\sigma_{\xb}^{-(\#{\omega})}|\Rv_{i+\omega}|^{1/2}\exp\acc{\frac12\mb_{i+\omega}\T\Rv_{i+\omega}\M\mb_{i+\omega}}\pth{\frac{\lambda}{1-\lambda}}^{\#\omega}
\eeq

Compared to a direct implementation using~\eqref{R}-\eqref{p}, several sources of computational saving can be found. The main one exploits the fact that both \Hv and $\Hv\T\Hv$ are Toeplitz matrix in the ``zero boundary'' case of finite convolution. As a consequence, neither $\Hv_{i+\omega}$, nor $\Hv_{i+\omega}\T\Hv_{i+\omega}$, nor $\Rv_{i+\omega}$ depend on the position $i$. A shorter notation where $i+\omega$ is replaced by $\omega$ will thus be adopted. More importantly, it becomes possible to precompute and store $\Rv_{\omega}\M$ and $|\Rv_{\omega}|$ for all values of $\omega$. An even more efficient scheme manipulates Cholesky factors as detailed in Table~\ref{algo:Ktuple}, such that $\Uv_{\omega}\T \Uv_{\omega} = \Sv_{\omega}$, $\Uv_{\omega}$ being an upper-triangular matrix.

\begin{table}[h]
\caption{Implementation of $K$-tuple sampling, Step~1(a),(b) in Table~\ref{Generalized Gibbs}.}
\label{algo:Ktuple}
\centering
\medskip
\bfmi[.65\linewidth]
\begin{algorithmic}
\FORALL{nonempty $\omega \subset \Omega=\{0,\ldots,K-1\}$} 
\STATE $L_{\omega} \gets \#\omega$
\STATE $\Sv_{\omega} \gets \sigma_{\epsilon}^{-2}\Hv_{\omega}\T\Hv_{\omega}+ \sigma_{\xb}^{-2} \texttt{eye}(L_{\omega})$ \COMMENT{$=\Rv_{\omega}\M$}
\STATE $\Uv_{\omega} \gets \texttt{chol}(\Sv_{\omega})$ \COMMENT{costs $\Oc(L_{\omega}^2)$ since $\Rv_{\omega}\M$ is Toeplitz}
\STATE $\alpha_{\omega} \gets |\Uv_{\omega}|\M$ \COMMENT{$=|\Rv_{\omega}|^{1/2}$}
\ENDFOR
\FOR{$i=1$ to $M-K+1$}
\STATE $\eb_{K,i} \gets \zb- \Hv_{-(i+\Omega)}\xb_{-(i+\Omega)}$
\FORALL{nonempty $\omega \subset \Omega=\{0,\ldots,K-1\}$} 
\STATE $c_{\omega} \gets \Uv_{\omega}\T\backslash \Hv_{\omega}\T\eb_{K,i}/\sigma_{\epsilon}^2$
\STATE $p_{i+\omega} \gets  \sigma_{\xb}^{-L_{\omega}} \alpha_{\omega} \exp(\norm{c_{\omega}}^2/2)(1/\lambda-1)^{-L_{\omega}}$ 
\ENDFOR
\STATE $p_{i+\emptyset} \gets 1$
\STATE \COMMENT{Step~1(a)}
\STATE sample $\qb_{i+\Omega}$ according to~\eqref{p}
\STATE \COMMENT{Step~1(b): $\xb_{i+\omega}\,\sim\,\mathcal{N}(\mb_{\omega},\Rv_{\omega})$}
\STATE $\xb_{i+\omega} \gets \Uv_{\omega}\backslash(c_{\omega} + \texttt{randn}(L_{\omega},1))$ 
\STATE $\xb_{i+\Omega\setminus\omega} \gets \zerob$ 
\ENDFOR
\end{algorithmic}
\efmi
\end{table}

Using an implementation of Table~\ref{algo:Ktuple}, the iteration numbers of Table~\ref{tab:K-tupleIter} can be converted in computer time. This confirms that despite the larger time taken per iteration, the K-tuple sampler may be faster for $K>1$ than for $K=1$. However, because of the exponentially increasing cost per iteration as a function of $K$, the best trade-off cannot be reached for large values of $K$. Here and after, we have limited our study to values of $K$ not larger than 4. Some critical results of the deconvolution problem of the Mendel's sequence in Fig.~\ref{fig:MCMCTimeshifted} are reported in Tab.~\ref{tab:K1-4}: including the overall computational costs using $4$ samplers in the $K$-tuple family and the corresponding number of iterations until convergence.
More complete simulation results are provided in Section~\ref{sec:Simulation results}. 
\btabl[htb]
\caption{Convergence rate of the K-tuple ($K=1,\ldots,4$) Gibbs scheme for the example in Fig.~\ref{fig:MCMCTimeshifted} measured in iteration numbers and time in seconds. Iteration numbers are noted in each case, e.g. $4600$ for the simple hybrid sampler and $900$ for the $2$-tuple , etc. \label{tab:K1-4}}
\centering
\btabu{|c||c|c|c|c|}
\hline
		&		$K =1$ & $K =2$ & $K =3$ & $K =4$ \\
\hline
Number of iterations & 	4600 & 900 & 700& 600 \\
\hline
Computation time (s) & 915 & 202 & 285 & 452 \\
\hline
\etabu
\etabl

\section{Partially marginal Gibbs sampler}

\label{sec:Partially marginal Gibbs sampler}
Here, another sampling scheme is proposed and studied, that indirectly tackles the same inefficiency by marginalizing out \xb in the step that samples \qb.

\subsection{Partial marginalization of \xb}
\label{subsec:Another Gibbs sampler}
In principle, one alternative to the approach in Eq.~\eqref{eqn:MonteCarlo} could be to generate a Markov chain on $\widetilde{\Theta}=(\qb,\hb,\sigma_{\epsilon}^2,\sigma_{h}^2, \lambda)$ marginalizing \xb, and then build the estimates in two steps: (1) the marginal maximum \apost estimate on each binary parameter $q_i$ and the posterior mean estimate on the continuous parameters $\widetilde{\Theta}\setminus \qb$, (2) the posterior mean estimate of \xb conditional on $(\widetilde{\Theta},\zb)$, that is a linear estimation problem. 
Furthermore, a Gibbs sampler with target distribution $P(\widetilde{\Theta}\I\zb)$ is likely to be more efficient than the hybrid sampler, particularly with respect to the sampling of \qb. Theoretical foundations are available in~\cite{Liu94b} and~\cite[Chapter 6.7]{Liu01} regarding the convergence rate of the so called \emph{collapsed Gibbs sampler}, with the marginalized equilibrium distribution $P(\widetilde{\Theta}\I\zb) = \int P(\Theta\I\zb)d\xb$. According to~\cite[Chapter 6.7]{Liu01}, the collapsed Gibbs sampler produces a substantial gain in terms of convergence rate if the marginalized parameters form highly dependent pairs with other parameters of the Markov chain~(see also~\cite{Cappe99} and references therein). This is exactly the case of the pair $\{\xb,\qb\}$ in the spike train deconvolution problem. Park \etal~\cite{Park09} also illustrates through three examples that such a sampler must be implemented with care to be sure that the desired stationary distribution is preserved. 

However, marginalizing \xb from the posterior distribution by $P(\widetilde{\Theta}\I\zb) = \int P(\Theta\I\zb)d\xb$ leads to practical difficulties.
Technically, the marginalization of \xb is only feasible on conditional laws of $(\qb, \lambda, \sigma_{\hb}^2)$. On the contrary, the conditional sampling of \hb becomes extremely difficult if \xb is integrated out, since $P(\hb\I\qb,\sigma_{\hb}^2,\sigma_{\epsilon}^2, \zb)$ is a multivariate, non Gaussian law with a complex structure. The same observation is made on the sampling of $\sigma_{\epsilon}^2$ if \xb is marginalized out.
Instead of a plain collapsed Gibbs sampler on $\widetilde{\Theta}$, the sampling scheme in Table~\ref{Di} is proposed to circumvent the direct conditional sampling of $\hb$ and $\sigma_{\epsilon}^2$.
Compared with the hybrid sampler, the main difference appears in Step~1, where \xb has been analytically integrated out. Subsection~\ref{subsec:Marginal posterior distribution} provides an efficient way to implement Step~1. 
\begin{table}[hbt]
\caption{The proposed partially marginal sampler.}
\label{Di}
\centering
\medskip
\bfmi[.5\linewidth]
\btabu{ll}
\textcircled{\tiny1} &(a) for each $i=1\ldots,M$,  \notag\\
&draw $q_i^{(k+1/2)}\I\qb_{1:i-1}^{(k+1/2)},\qb_{i+1:M}^{(k)},\hb^{(k)},\sigma_{\epsilon}^{(k)},\lambda^{(k)}, \zb$\\ 
                     &(b) draw $\xb^{(k+1/2)} \I\qb^{(k+1/2)},\hb^{(k)},\sigma_{\epsilon}^{(k)}, \zb$\\
\textcircled{\tiny2}~&(a) draw $\yb^{(k+1)}$ according to Table~\ref{algo:TimeShift}\\
                     &(b) draw $\hb^{(k+1)} \I\xb^{(k+1)},\sigma_{\hb}^{(k)},\sigma_{\epsilon}^{(k)}, \zb$ \\
\textcircled{\tiny3} &draw $\sigma_{\epsilon}^{(k+1)}\I\xb^{(k+1)},\hb^{(k+1)}, \zb$\\ 
\textcircled{\tiny4} &draw $\lambda^{(k+1)}\I\qb^{(k+1)}$\\
\textcircled{\tiny5} &draw $\sigma_{\hb}^{(k+1)}\I\hb^{(k+1)}$ 
\etabu
\efmi
\end{table}

We give here a mathematical justification that the proposed partially marginal sampler identifies with a Gibbs sampler of the joint posterior law (Eq.~\eqref{eqn:JointPostdistribution}, that includes \xb), for a particular scanning scheme. The techniques are similar to those denoted as \emph{marginalization} and \emph{trimming} in~\cite{Park09}.

Let $(q_1,\xb), (q_2,\xb),\ldots,(q_M,\xb), \hb,\sigma_{\epsilon}^2,\lambda,\sigma_{\hb}^2$ be the scanning sequence of the Gibbs sampler. For the sampling of couples $(q_i,\xb)$, a two-stage procedure is considered by applying the Bayes rule $P(\xb, q_i \I \Theta \setminus \{\xb, q_i\} ,\zb) = P(q_i \I \Theta \setminus \{\xb, q_i\} ,\zb)P(\xb \I \Theta \setminus \xb, \zb)$: 
\bit
\item draw $q_i$ according to $P(q_i \I \qb_{1:i-1},\qb_{i+1:M},\hb,\sigma_{\epsilon}^2,\lambda,\zb)$;
\item draw \xb according to $P(\xb \I \qb,\hb,\sigma_{\epsilon}^2,\lambda,\zb)$.
\eit

The second stage is therefore repeated $M$ times within each iteration of the sampler, while all but the last one are useless since the correspondingly sampled values of \xb do not enter any subsequent operation. The $M-1$ corresponding sampling operations can thus be skipped and the resulting sequence of operations coincides exactly with the partially marginal sampler of Table~\ref{Di}. 
Another mathematical proof is given in~\cite{Ge08} by considering Table~\ref{Di} as a collapsed Gibbs sampler with  target distribution $P(\widetilde{\Theta}\I\zb)$ and then checking the \emph{invariant condition} of the related Markov chain. 

\begin{figure}[htb]
  \centering
\btabu{@{}c@{~~}c@{}}
  \figc[width=.47\linewidth]{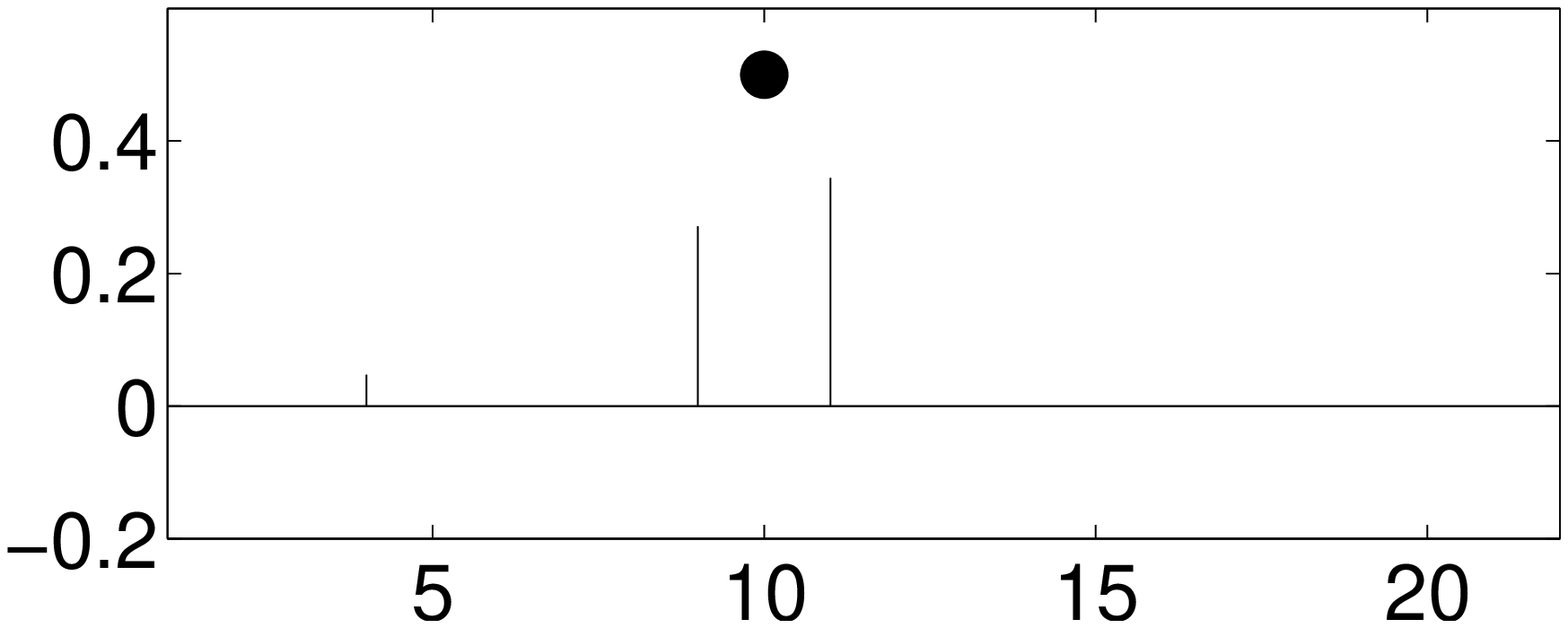} & \figc[width=.47\linewidth]{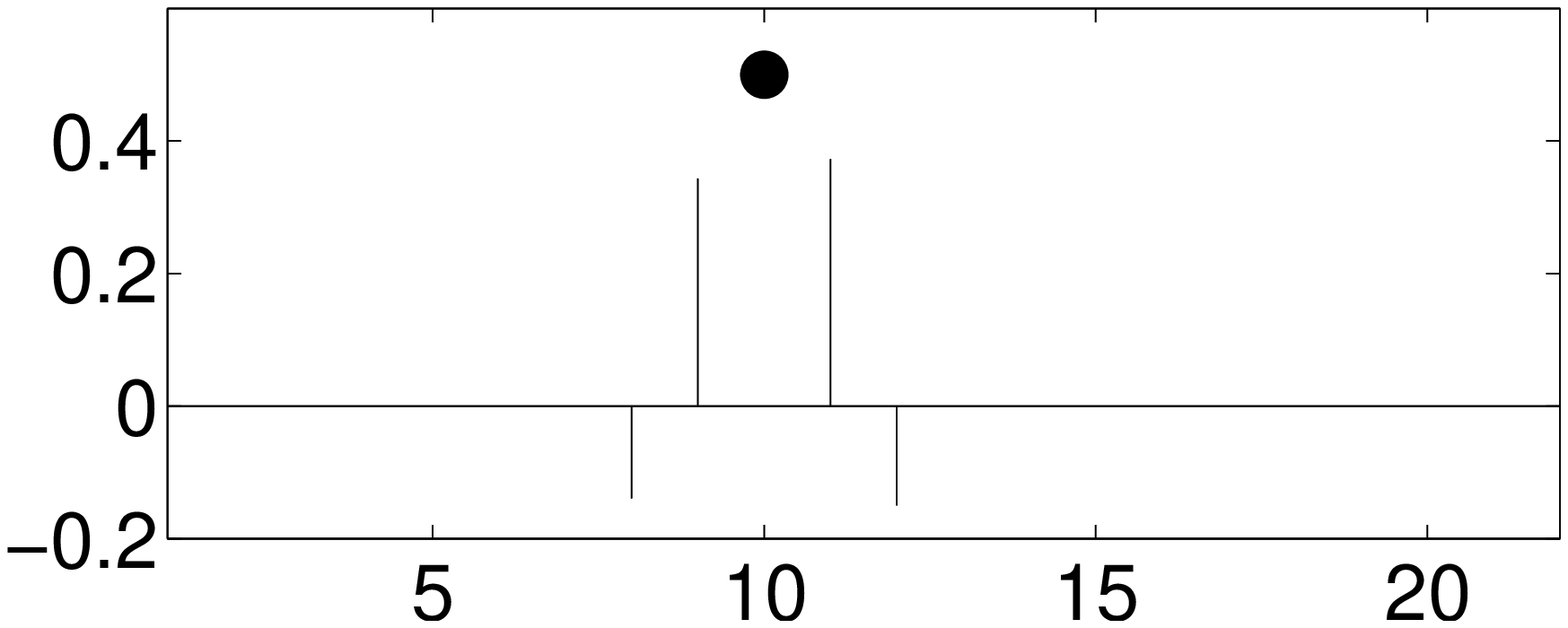} \\
{\footnotesize  (a) 1st iteration} & {\footnotesize (b) 18th iteration}\\[2mm]
  \figc[width=.47\linewidth]{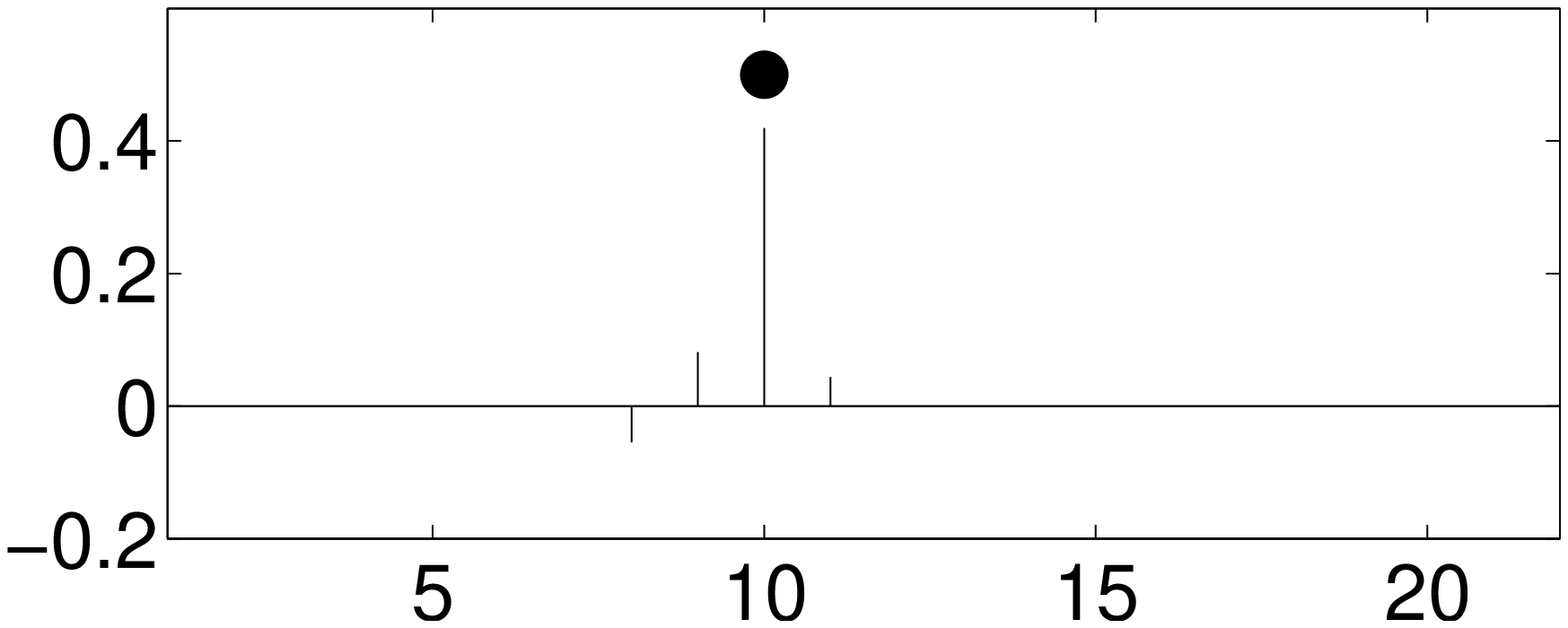} & \figc[width=.47\linewidth]{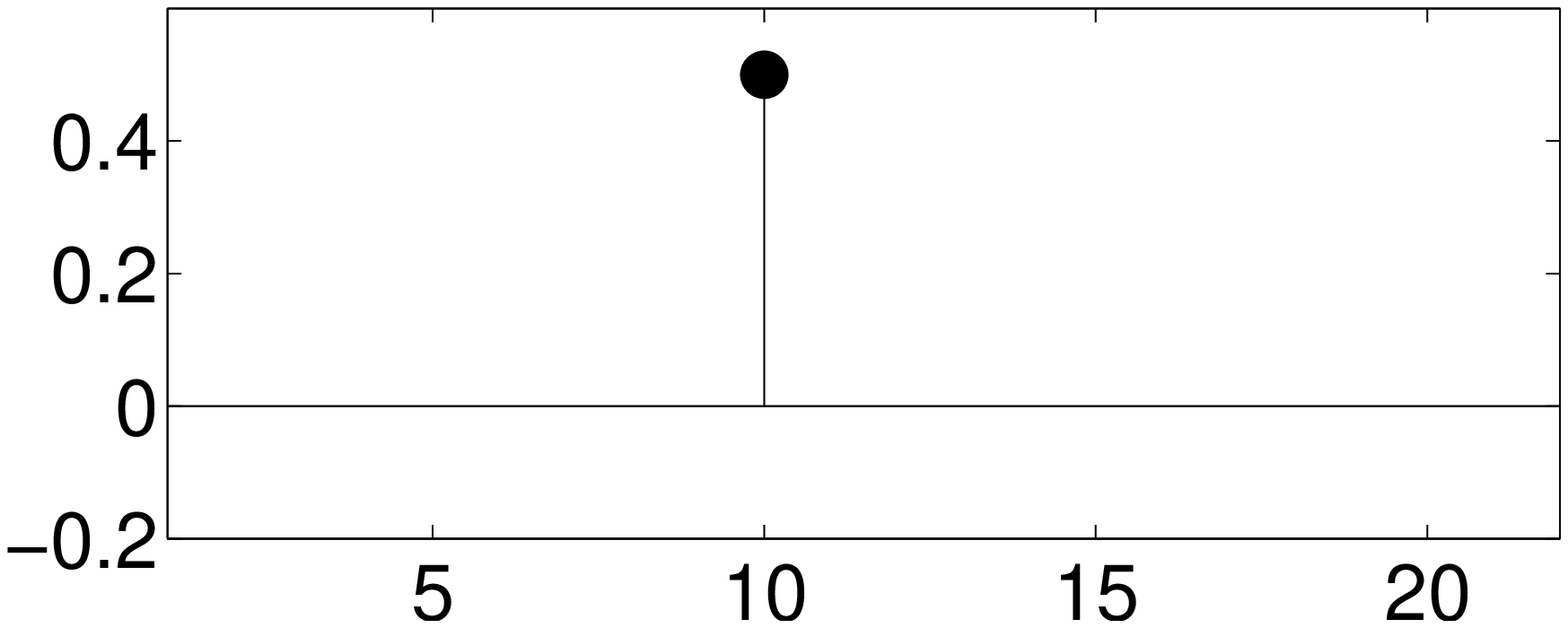} \\
{\footnotesize  (c)  19th iteration} & {\footnotesize (d) 20th iteration}\\
\etabu
\caption{Sampled BG sequences obtained by our partially marginal scheme in the example of Figure~\ref{fig:Configs}.
The Markov chain escapes rapidly the initial configuration neighborhood.
}
\label{fig:Escape}
\end{figure}

In the same conditions than that of Table~\ref{tab:K-tupleIter}, the partially marginal sampler takes about $20$ iterations.
Figure~\ref{fig:Escape} illustrates the way it escapes from a local optimal configuration within an acceptable number of iterations in the example of Figure~\ref{fig:Configs}. The true configuration is reached after 20 iterations. Moreover, it is observed from the configurations obtained at the 18th and 19th iterations that our scheme is able to radically modify \xb in one single step, a characteristic directly related to the marginalization of \xb in Step~1(a) of the sampler.

\subsection{Marginal posterior distribution}
\label{subsec:Marginal posterior distribution}
In this subsection, a numerical implementation of Step~1 of Table~\ref{Di} is proposed, inspiring from
the recursive method in~\cite{Champagnat96a} to evaluate $P(q_i=0,1|\widetilde{\Theta}\setminus q_i;\zb)$ sequentially. While the latter is based on the storage and update of an $L \times L$ ($L=\sum_i q_i$) matrix, we introduce an even less burdensome strategy by handling the Cholesky factor matrix. The issue of the complexity per iteration is dealt here, whereas the overall computational costs are compared in Section~\ref{sec:Simulation results} using simulation tests.
 
First, the conditional posterior distribution of $q_i$ takes the following form, as detailed in~\cite{Champagnat96a}:
\balx
P(q_i|\widetilde{\Theta}\setminus q_i; \zb)
\propto |\Bv|^{-1/2}\exp\left(-\frac12\zb\T\Bv\M\zb\right)\pth{\frac{\lambda}{1-\lambda}}^{q_i}\nonumber
\propto \exp\left(-f(q_i)/2\right)
\ealx
where
\bal
\Bv &= \Hv\diag(\qb)\Hv\T + \sigma_{\epsilon}^2\Iv_N, \label{eqn:B}\\
f(q_i) &= \zb\T\Bv\M\zb + \log\bars{\Bv} + 2q_i\log(1/\lambda-1). \label{eqn:f}
\eal
After normalization, the marginal probability of $q_i$ reads:
\beqx
P(q_i\I\widetilde{\Theta}\setminus q_i, \zb)=\pth{1+\exp\left(-
(f(1-q_i)-f(q_i))/2\right)}\M,
\eeqx
which reduces to the evaluation of $f(1-q_i)-f(q_i)$, using the following steps~\cite{Champagnat96a}:
\bal
\tau_i &= \delta_i + \sigma_{\epsilon}^{-2} \hb_i\T\widetilde{\Bv}_0\M \hb_i \\
\widetilde{\Bv}_i\M &= \widetilde{\Bv}_0\M - \sigma_{\epsilon}^{-2} \widetilde{\Bv}_0\M\hb_i \tau_i\M \hb_i\T \widetilde{\Bv}_0\M \label{eqn:deltaL}\\
\stdbars{\widetilde{\Bv}_i\M} &= \stdbars{\widetilde{\Bv}_0\M}\delta_i\tau_i \label{eq:DetBi}
\eal
with $\widetilde{\Bv} = \Bv/\sigma_{\epsilon}^2$. $\delta_i=\pm1$ depends on whether $1$ is added or removed at $q_i$ and $\widetilde{\Bv}_i$ and $\widetilde{\Bv}_0$ differ only at $q_i$. Let us note that $\widetilde{\Bv}_0, \widetilde{\Bv}_i > 0$, so that it can be deduced from~\eqref{eq:DetBi} that $\tau_i$ has the same sign as $\delta_i$.

Further simplifications are also introduced in~\cite{Champagnat96a} by exploiting the sparse nature of $\widetilde{\Bv}$ and applying the matrix inversion lemma. Noticing that the rank of $\Hv\diag(\qb)\Hv\T$ is only $L$, $\widetilde{\Bv}$ takes the alternate form $\widetilde{\Bv} = \sigma_{\epsilon}^{-2}\Gv\Gv\T+\Iv$, where $\Gv=\Hv\Dv$ is full rank, and $\Dv$ made of the nonzero columns of $\diag(\qb)$. By applying the matrix inversion lemma, we have
\bal
\widetilde{\Bv}\M = \Iv - \sigma_{\epsilon}^{-2} \Gv\Cv\M\Gv\T \label{eqn:Binv}
\eal
where $\Cv = \sigma_{\epsilon}^{-2}\Gv\T\Gv+\Iv$ is a $L\times L$ matrix.
Therefore, $\Cv\M$ ($L\times L$) can be stored and updated instead of $\widetilde{\Bv}\M$ ($N\times N$). 
In the case of adding a spike at $q_i$, the formula for the update of $\Cv_i\M$ up to a matrix permutation operation is the following: 
\bal
\bb &= -\sigma_{\epsilon}^{-2} \tau_i\M \Cv_0\M\Gv_0\T\hb_i \label{eqn:updateb} \\
\Cv_i\M &= \left[
\barr{cc}
\Cv_0\M+\bb\tau_i\bb\T & \bb \\[3pt]
\bb\T & \tau_i\M
\earr\label{eqn:updateC}
\right]
\eal
The case of removing a spike is straightforward since we have the following:
\bal \label{eq:UpdateCminus}
\Cv_0\M &= \left[
\barr{cc}
\Cv_i\M+\bb\tau_i\bb\T & \bb \\[3pt]
\bb\T & \tau_i\M
\earr
\right]
\eal
Notice that in the case of spike removal, \bb and $\tau_i$ are extracted from $\Cv_0\M$ rather than calculated as in Eq.~\eqref{eqn:updateb}.  
These operations are invariant up to a matrix permutation in cases of adding/removing a spike at an arbitrary location.  

Finally, we propose to further reduce the computation and memory load, by updating and storing the Cholesky factor $\Fv$ of $\Cv\M$ ($\Fv\T\Fv = \Cv\M$ for an upper-triangular matrix \Fv). 
It is derived from Equations~\eqref{eqn:deltaL} and~\eqref{eqn:Binv} that:
\balx
\tau_i&=\delta_i+\sigma_{\epsilon}^{-2}\norm{\hb_i}^2-\sigma_{\epsilon}^{-4}\stdnorm{\Fv_0\Gv_0\T\hb_i}^2,\\
\zb\T\widetilde{\Bv}_0\M\hb_i&=\zb\T\hb_i-\sigma_{\epsilon}^{-2}(\Fv_0\Gv_0\T\zb)\T(\Fv_0\Gv_0\T\hb_i)
\ealx
and from Eq.~\eqref{eqn:f} and~\eqref{eqn:deltaL} that:
\bal\label{eq:deltaF}
f(1-q_i)-f(q_i)= \log(\delta_i\tau_i)-\sigma_{\epsilon}^{-4}\tau_i\M(\zb\T\widetilde{\Bv}_0\M\hb_i)^2+2\delta_i\log\pth{\UNsur{\lambda}-1}
\eal

Two extra advantages of the Cholesky decomposition are:
\bit  
\item since the conditional law in Step~1(b) reads 
\beq\label{eq:UpdateX}
(\xb\I\wt{\Thetab},\zb)\, \sim \,\Nc(\sigma_{\epsilon}^{-2}\Cv\M\Gv_0\T\zb,\Cv\M),
\eeq 
sampling \xb costs $\Oc(L^2)$ using \Fv (See Table~\ref{algo:Marginal}); 
\item the incremental update of $\Fv$ also costs $\Oc(L^2)$ using a {\it rank-$1$ Cholesky update} method. 
The examples of adding and removing a spike at the last location is given in the following for notational simplicity. It is however detailed in the appendix that unlike adding a spike at at an arbitrary location, the case of removing a spike at an arbitrary location requires two operations of Cholesky rank-$1$ update. 

For adding a spike at the last location, $\tau_i>0$ and equation~\eqref{eqn:updateC} also reads:
\beq\label{eq:CholUpdateC}
\Cv_i\M = \cro{
\barr{c@{\qquad}c}
\Fv_0\T & \zerob \\
\zerob\T & 0
\earr
}
\left[
\barr{c@{\qquad}c}
\Fv_0 & \zerob \\
\zerob\T & 0
\earr
\right]+\left[
\barr{c}
\bb\tau_i^{1/2} \\
\tau_i^{-1/2} 
\earr
\right]
\left[
\barr{c}
\bb\tau_i^{1/2} \\
\tau_i^{-1/2}
\earr
\right]\T
\eeq
\eit 
In the case of removing a spike at the last location, from Eq.~\eqref{eq:UpdateCminus} and the definition of Cholesky decomposition, we have:
\bal
\Fv_0\T \Fv_0 &= \left[
\barr{cc}
\Fv_i\T\Fv_i+\bb\tau_i\bb\T & \bb \\
\bb\T & \tau_i\M
\earr
\right] \\ 
\Fv_0 &= \left[
\barr{cc}
\Sv_0 &  \ab \\
\zerob\T  &  m
\earr
\right]
\eal
where $\Sv_0 $ is also an upper triangular matrix. It is then directly deduced that: 
\bal\label{eqn:UpdateFb}
\Sv_0\T \ab &= \bb \\
\ab\T \ab + m^2 &= \tau_i\M \\
\Sv_0\T \Sv_0 &= \Fv_i\T\Fv_i +\bb\tau_i\bb\T	 
\eal
Let $\texttt{cholupdate}(\Av, \db, '\pm')$ be the method that updates in $\Oc(L^2)$ the Cholesky factor of $\Av\T\Av\pm \db\db\T$, where \Av is already an upper triangular matrix. Then, $\Fv_i = \texttt{cholupdate}(\Sv_0, \bb\tau_i^{1/2}, '-')$. 
\begin{table}[htb]
\caption{Step~1 of the partially marginal sampling scheme. The update of \Fv is given here up to a matrix permutation.}
\label{algo:Marginal}
\centering
\medskip
\bfmi[.9\linewidth]
\begin{algorithmic}
\STATE Initialize \Fv and \Gv
\STATE \COMMENT{Step~1(a): sequential sampling of $q_i$}
\FOR{$i=1$ to $M$}
\STATE $\delta_i \gets (-1)^{q_i}$
\STATE $\tau_i \gets \delta_i+\sigma_{\epsilon}^{-2}\norm{\hb}^2-\sigma_{\epsilon}^{-4}\stdnorm{\Fv\Gv\T\hb_i}^2$
\STATE $\phi \gets \zb\T\hb_i-\sigma_{\epsilon}^{-2}(\Fv\Gv\T\zb)\T(\Fv\Gv\T\hb_i)$ \quad \quad \COMMENT{$\phi=\zb\T\Bv\hb$}
\STATE $\Delta f \gets \texttt{log}(\delta_i\tau_i)-\sigma_{\epsilon}^{-4}\tau_i\M\phi^2+2\delta_i\log\pth{1/{\lambda}-1}$ \quad \COMMENT{see Eq.\eqref{eq:deltaF}}
\STATE sample $u \sim \Uc([0,1])$ \COMMENT{Uniform distribution in $[0,1]$}
\IF{$u > (1+\exp(-\Delta f/2))\M$} 
\STATE $q_i \gets q_i + \delta_i$ \quad \COMMENT{$q_i + \delta_i$ is accepted}
\IF{$\delta_i=1$}
\STATE $\bb \gets -\sigma_{\epsilon}^{-2} \tau_i\M \Fv\T\Fv\Gv\T\hb_i$ \quad \COMMENT{see Eq.~\eqref{eqn:updateb}}
\STATE $\Fv \gets \texttt{cholupdate}([\Fv\ \zerob; \zerob\T\ 0],[\bb \tau_i^{1/2};\tau_i^{-1/2}], '+')$ \quad \COMMENT{see Eq.\eqref{eq:CholUpdateC}}
\ELSE[the case of $\delta_i= -1$ removing a spike, see appendix for detail]
\STATE $L \gets \text{sum}(\qb)$
\STATE $\bb \gets \Fv\T\Fv(:,i), \vb \gets \bb_{-i}, \tau\M \gets \bb_i$
\IF[other than the last location, an extra update needed]{$i< L$}
\STATE $\eb \gets \Fv(i, 1+i:L)\T$
\STATE $\Fv(1+i:L, 1+i:L) \gets \texttt{cholupdate}(\Fv(1+i:L, 1+i:L), \eb, '+')$ \quad \COMMENT{$\Oc((L-i)^2)$}
\ENDIF
\STATE $\Fv(i, :) \gets [], \Fv(:, i) \gets []$ \COMMENT{Eliminate the $i$th row and column}
\STATE $\Fv \gets \texttt{cholupdate}(\Fv, \tau^{1/2}\vb, '-')$ \quad \COMMENT{$\Oc(L^2)$}
\ENDIF
\STATE $\Gv \gets \Hv(:\,,\texttt{find}(\qb))$
\ENDIF \quad \COMMENT{No change of $(q_i,\Fv)$ otherwise}
\ENDFOR
\STATE \COMMENT{Step~1(b): sampling of \xb}
\STATE $\xb \gets \zerob$ 
\STATE $L \gets \text{sum}(\qb)$
\STATE $\xb(\texttt{find}(\qb)) \gets \Fv\T\left(\sigma_{\epsilon}^{-2}\Fv\Gv\T\zb + \texttt{randn}(L,1)\right)$ \quad \COMMENT{see Eq.~\eqref{eq:UpdateX}}
\end{algorithmic}
\efmi
\end{table}

\section{Simulation results}

\label{sec:Simulation results}

\subsection{Convergence diagnostic}
In order to compare empirically the convergence speed of the different samplers, we have resorted to Brooks and Gelman's iterated graphical method to assess convergence ~\cite{Brooks98}. This diagnostic method is based upon the covariance estimation of $m$ independent Markov chains $\{\Phi_{jt}, j=1,\ldots,m;t=1,\ldots,n\}$ of equal length $n$. Let $\overline{\Phi}_{j.}$ (respectively, $\overline{\Phi}_{..}$) denote the local (respectively, global) mean of the chains. The intra-chain and inter-chain variances are defined as covariance matrix averages:
\balx
\Vv_{\text{intra}} &= \frac1{m(n-1)}\sum_{j=1}^m\sum_{t=1}^n(\Phi_{jt}-\overline{\Phi}_{j.})(\Phi_{jt}-\overline{\Phi}_{j.})\T \nonumber\\
\Vv_{\text{inter}} &= \frac1{m-1}\sum_{j=1}^m(\overline{\Phi}_{j.}-\overline{\Phi}_{..})(\overline{\Phi}_{j.}-\overline{\Phi}_{..})\T
\ealx
that characterize the convergence behavior.
Brooks and Gelman~\cite{Brooks98} proposed to evaluate the \emph{multivariate potential scale reduction factor} (MPSRF):
$$
\text{MPSRF}=\frac{n-1}{n}+\frac{m+1}{m}\lambda(\Vv_{\text{intra}}\M\Vv_{\text{inter}}) 
$$
where $\lambda(\cdot)$ returns the largest eigenvalue of the covariance matrix. Convergence is diagnosed when MPSRF is close to one (\eg $\text{MPSRF}<1.2$ as proposed in~\cite{Brooks98}). In order to get a graphical evolution of the convergence, each chain is divided into batches of $b$ samples, and the MPSRF is calculated upon the second halves of the Markov chains $\{\Phi_{jt}\}, t=1,\ldots, kb$ of increasing lengths $kb$, while the first halves are regarded as a burn-in period. From an empirical basis, it is also advised to select $b\approx n/20$ in~\cite{Brooks98}.  

\subsection{Simulation tests}
\begin{figure}[htb]
\mongraphe{\includegraphics[width=160mm]{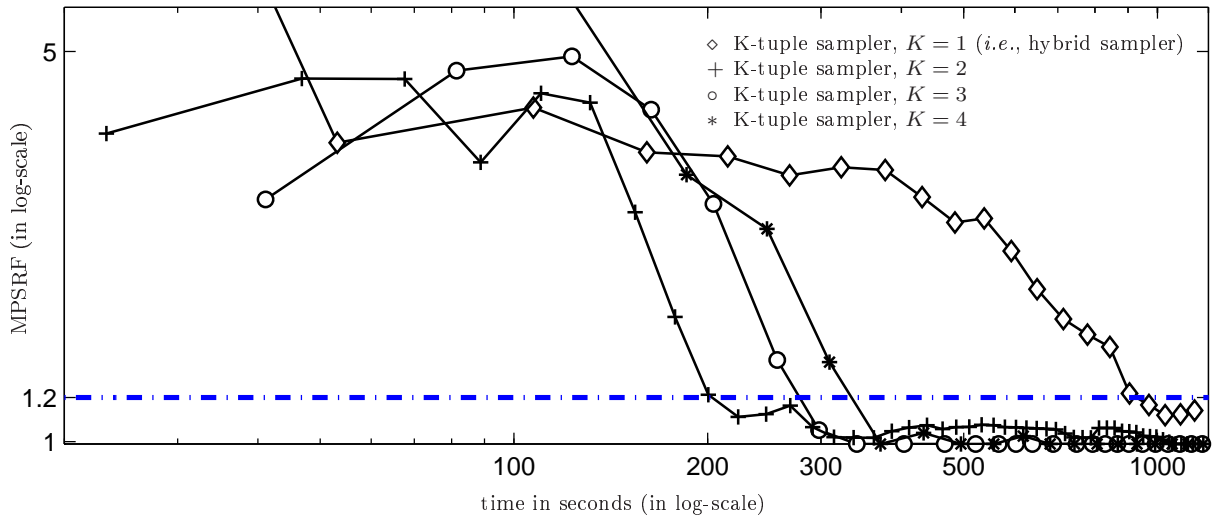}}%
{\hspace{60mm}
	{\footnotesize time in seconds (in log-scale)}
}%
{\hspace{22mm}
	{\footnotesize MPSRF (in log-scale)}
}%
{\hspace{90mm}\vspace{50mm}%
	\renewcommand{\arraystretch}{0.8}
	\begin{tabular}{@{}l@{~}l@{}}
		{\scriptsize $\diamond$}      &{\scriptsize K-tuple sampler, $K=1$ (\ie hybrid sampler)}\\
		{\scriptsize $+$}             &{\scriptsize K-tuple sampler, $K=2$}\\
		{\scriptsize $\circ$}         &{\scriptsize K-tuple sampler, $K=3$} \\
		{\scriptsize $\ast$}          &{\scriptsize K-tuple sampler, $K=4$}
	\end{tabular}
}
\caption{Evolution of MPSRF for the $K$-tuple sampler family ($K=1,2,3$) in the case of Mendel's sequence ($\mathrm{SNR}=12.80$\,dB).\label{fig:12dBCompare}}
\end{figure}

A test scenario is designed here to compare the generalized K-tuple Gibbs sampler (including the classical hybrid sampler in the case $K=1$) of Section~\ref{sec:Generalized Gibbs on K-tuple variables} and the partially marginal sampler of Section~\ref{sec:Partially marginal Gibbs sampler}, in terms of robustness with respect to different random initial conditions. Time-shift and scale ambiguities are taken into account in all the methods, as described in Section~\ref{subsec:Classical Gibbs Sampling}.
Let us reexamine the example of the benchmark sequence (Mendel's sequence, for which $M=300$ and $\lambda=0.1$.) in Figure~\ref{fig:MCMCTimeshifted}.
To concentrate on the convergence quality of the Bernoulli sequence, we evaluate the MPSRF evolutions of $\{\qb^{(k)}\}$ for ten independent Markov chains, \ie $m=10$. 
Figure~\ref{fig:12dBCompare} is compares the convergence rate of $K$-tuple samplers for different values of $K$: while MPSRF falls under $1.2$ after $515$ seconds of simulation for the classical hybrid sampler ($K=1$), this threshold is reached after about $285$ seconds for the 3-tuple sampler and in less than 202 seconds for the 2-tuple sampler.
It is interesting to observe that both the $3$ and $4$-tuple sampler have a more stable convergence behavior in the convergence zone while the hybrid sampler presents the most oscillating MPSRF values. These results confirm the discussion in Section~\ref{sec:Generalized Gibbs on K-tuple variables}: augmenting $K$ improves the convergence quality in terms of iteration numbers 
while in the meantime increases the computational load per iteration. The overall performance relies on the trade-off between the two criteria. Thus, the best compromise in the $K$-tuple sampler {\it series} is reached for $K=2$ in the given example.  

\begin{figure}[htb]
\mongraphe{\includegraphics[width=160mm]{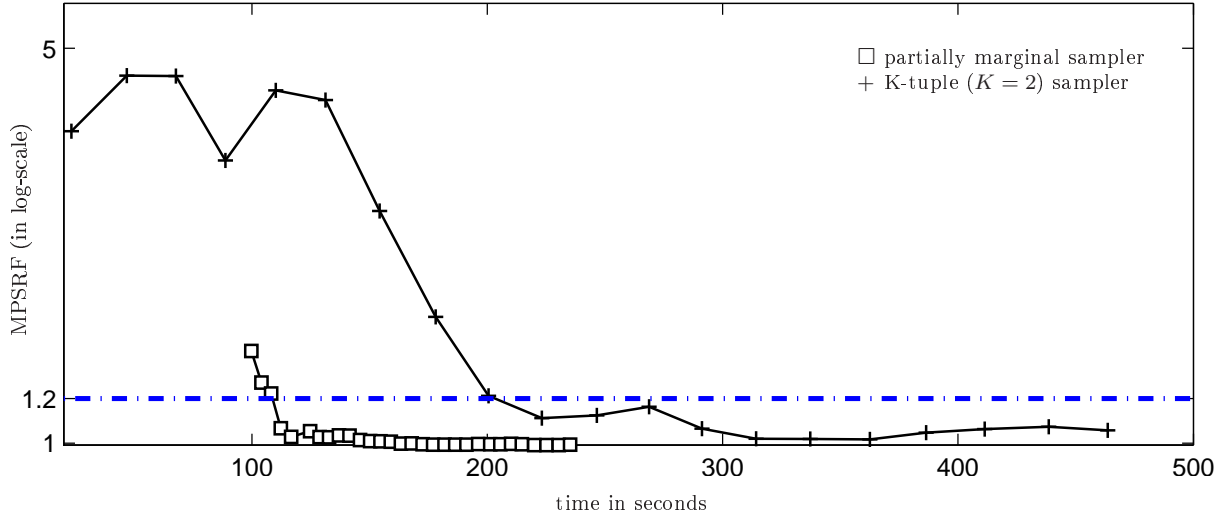}}%
{\hspace{70mm}
	{\footnotesize time in seconds}
}%
{\hspace{22mm}
	{\footnotesize MPSRF (in log-scale)}
}%
{\hspace{110mm}\vspace{55mm}%
	\renewcommand{\arraystretch}{0.8}
	\begin{tabular}{@{}l@{~}l@{}}
		{\scriptsize $\square$}      &{\scriptsize partially marginal sampler}\\
		{\scriptsize $+$}             &{\scriptsize K-tuple ($K=2$) sampler}
	\end{tabular}
}
\caption{Evolution of MPSRF for the partially marginal sampler and the $2$-tuple sampler in the case of Mendel's sequence. MPSRF has been computed every 100 samples for both methods.  \label{fig:GibbsTypeVS2Tuple}}
\end{figure}

Next we compare the convergence diagnostic results on the $2$-tuple sampler and the partially marginal sampler. Each MPSRF in Figure~\ref{fig:GibbsTypeVS2Tuple} is evaluated for every $b=100$ iterations. Altogether 300 iterations (116 seconds of simulation) are required for the partially marginalized sampler to reach convergence, not to mention the overall lower level of MPSRF in the convergence zone compared with the $2$-tuple sampler.   
Thus the partially marginalized sampler takes fewer iterations and less time to converge than the best selected $K$-tuple sampler in Figure~\ref{fig:12dBCompare}. 
We also point out that the asymptotic complexity per iteration of the partially marginal sampler is lower than that of the $2$-tuple sampler in the given simulation example, though the contrary is observed for the first 100 samples (corresponding to the heating period) in which the number of detected spikes $L$ is relatively important. 

\begin{figure}[htb]
  \centering
\setlength{\tabcolsep}{0pt}
\btabu{cc}
   \figc[height=5cm]{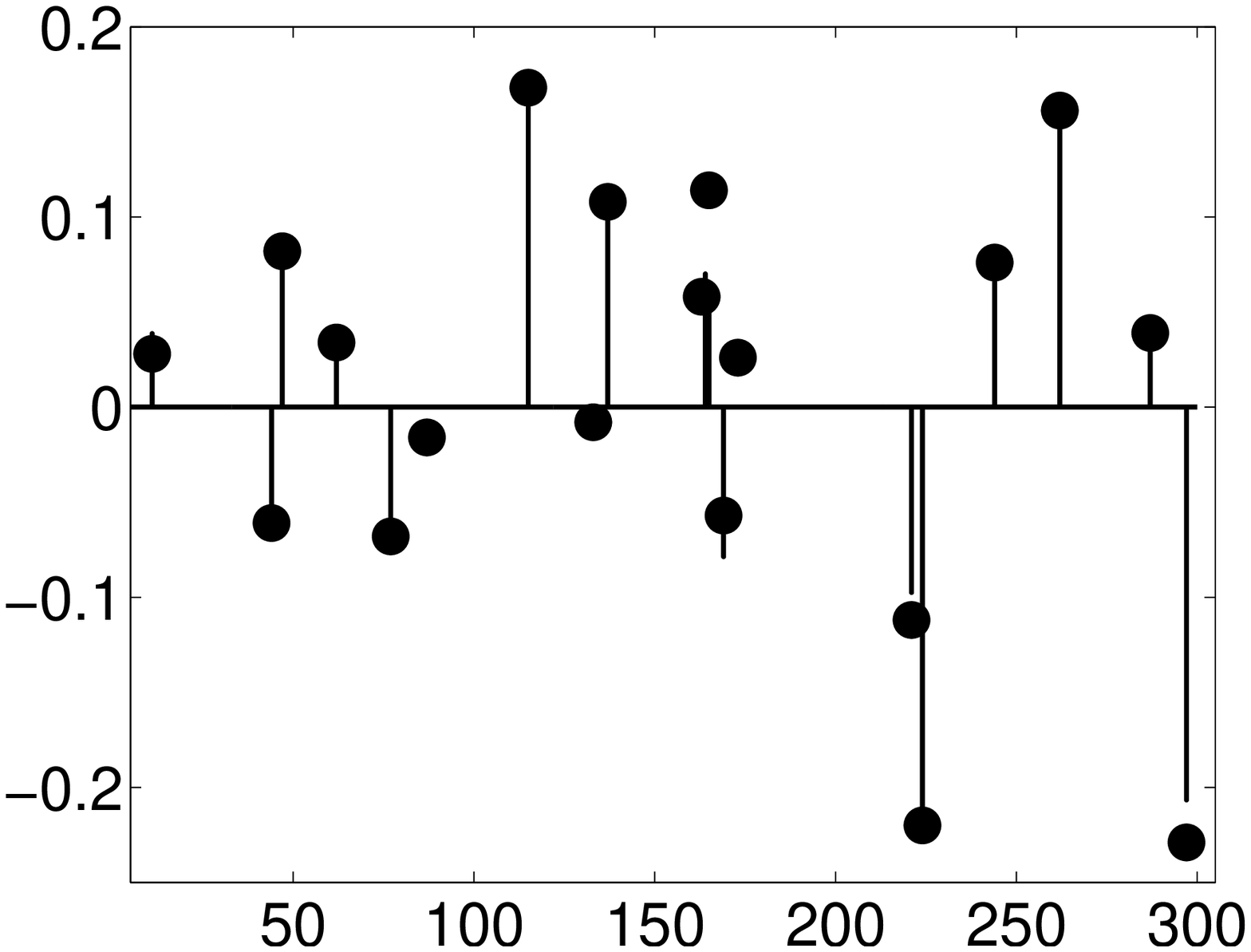}
&  \figc[height=5cm]{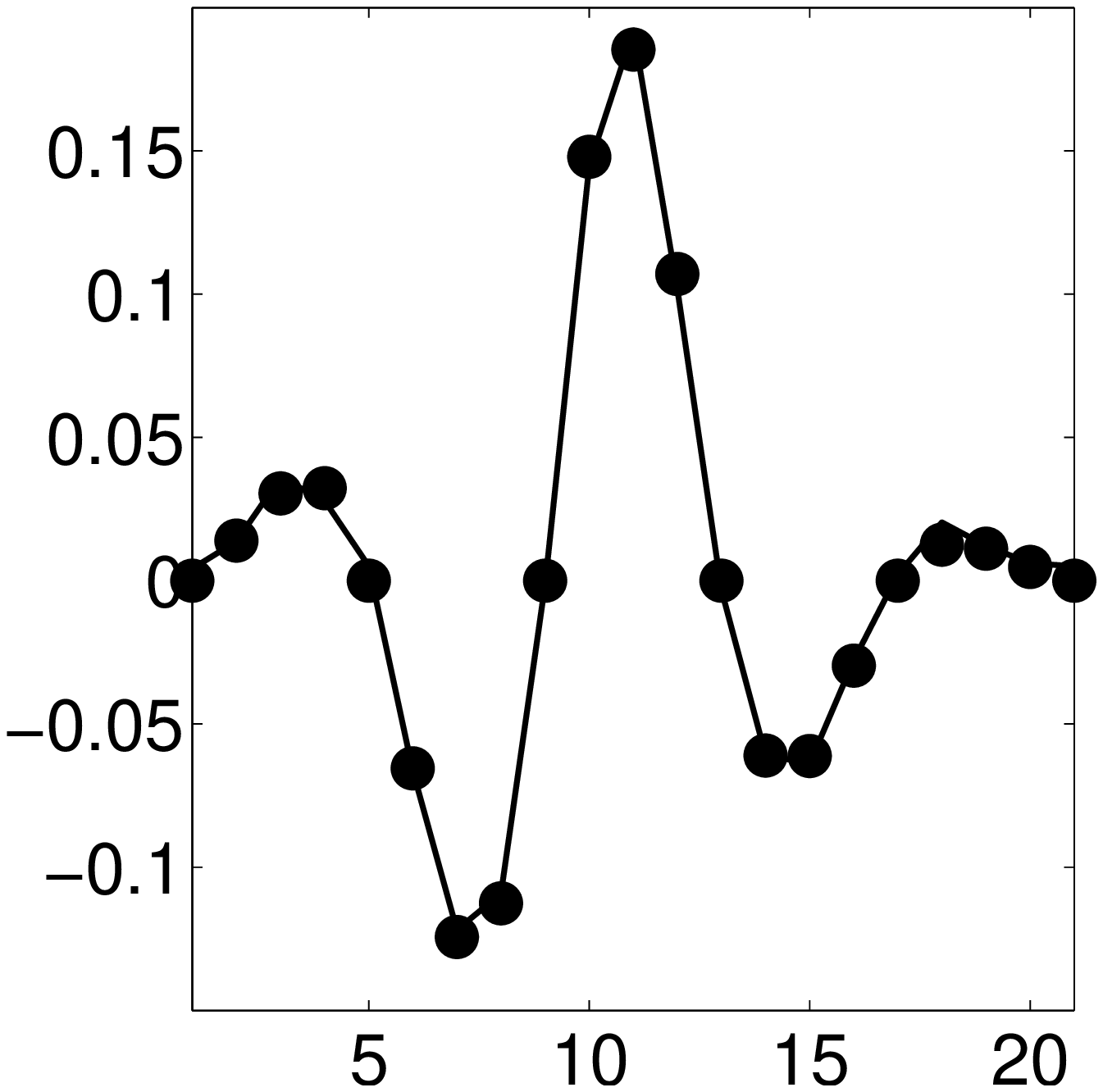} \\
		{\footnotesize estimation of \xb}
&  	{\footnotesize estimation of \hb}
\etabu
\caption{Estimation result using the proposed sampler, under the same conditions as in Figure~\ref{fig:MCMCTimeshifted}. The 10 independent Markov chains yield undistinguishable estimates.}
\label{fig:MCMCSLMR}
\end{figure}

To further illustrate its robustness, deconvolution results obtained by the partially marginal sampler are shown in Figure~\ref{fig:MCMCSLMR} under identical initial conditions as in Figure~\ref{fig:MCMCTimeshifted}. Only one of the results is reported, as they are undistinguishable from each other (perfectly consistent estimations). And this is true for all the samplers once their MPSRF values drop below the $1.2$ threshold. 


\begin{figure}[htb]
\mongraphe{\hspace{5mm} \includegraphics[height=90mm]{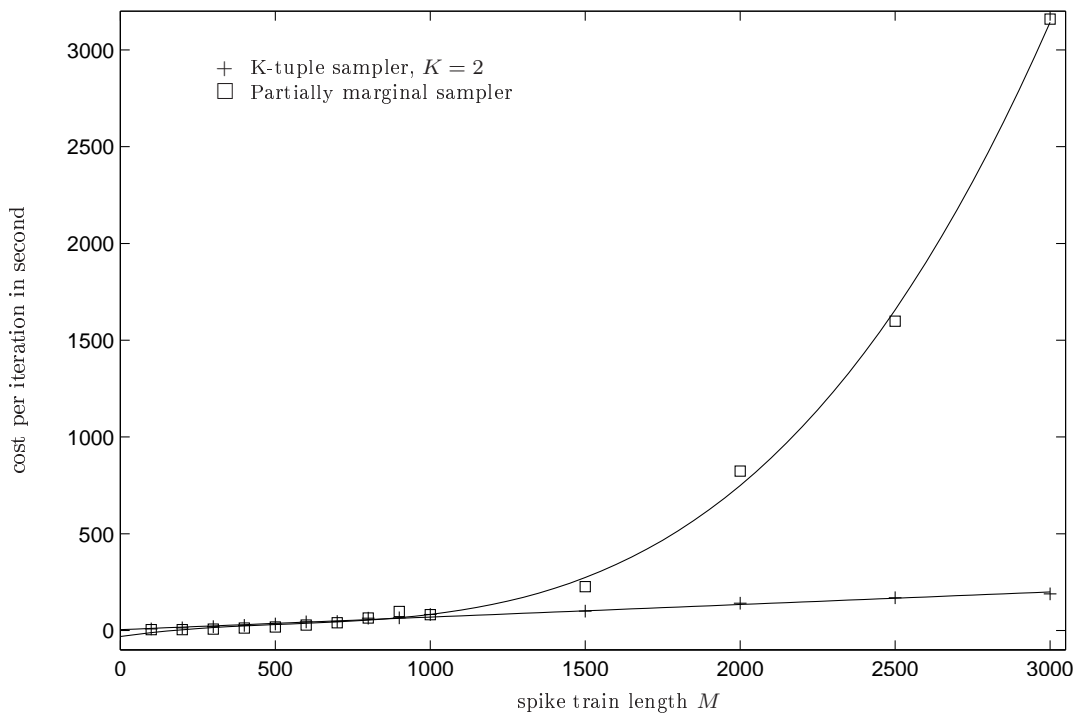}}%
{\hspace{65mm} {\footnotesize spike train length $M$}}
{\hspace{30mm} {\footnotesize cost per iteration in second}}
{\hspace{25mm}\vspace{80mm}
	\renewcommand{\arraystretch}{0.8}
	\begin{tabular}{@{}l@{~}l@{}}
		{\scriptsize $+$ }              &{\scriptsize K-tuple sampler, $K=2$}\\
		{\scriptsize $\square$}         &{\scriptsize Partially marginal sampler}
	\end{tabular}
}
\caption{Computation cost per iteration comparison between the $2$-tuple Gibbs sampler and the partially marginal sampler. The two methods are tested on the same simulated observation data \zb, for which the length of the spike train ranges from $100$ to $3000$ and the Bernoulli parameter $\lambda$ is fixed to $10\%$. Polynomial interpolations (of degree $1$ and $2$ respectively) are traced to show the quasi-linear and -quadratic evolution of complexity per iteration with respect to the spike train length $M$. }
\label{fig:CostIter}
\end{figure}

\begin{figure}[htb]
\mongraphe{\hspace{5mm} \includegraphics[height=100mm]{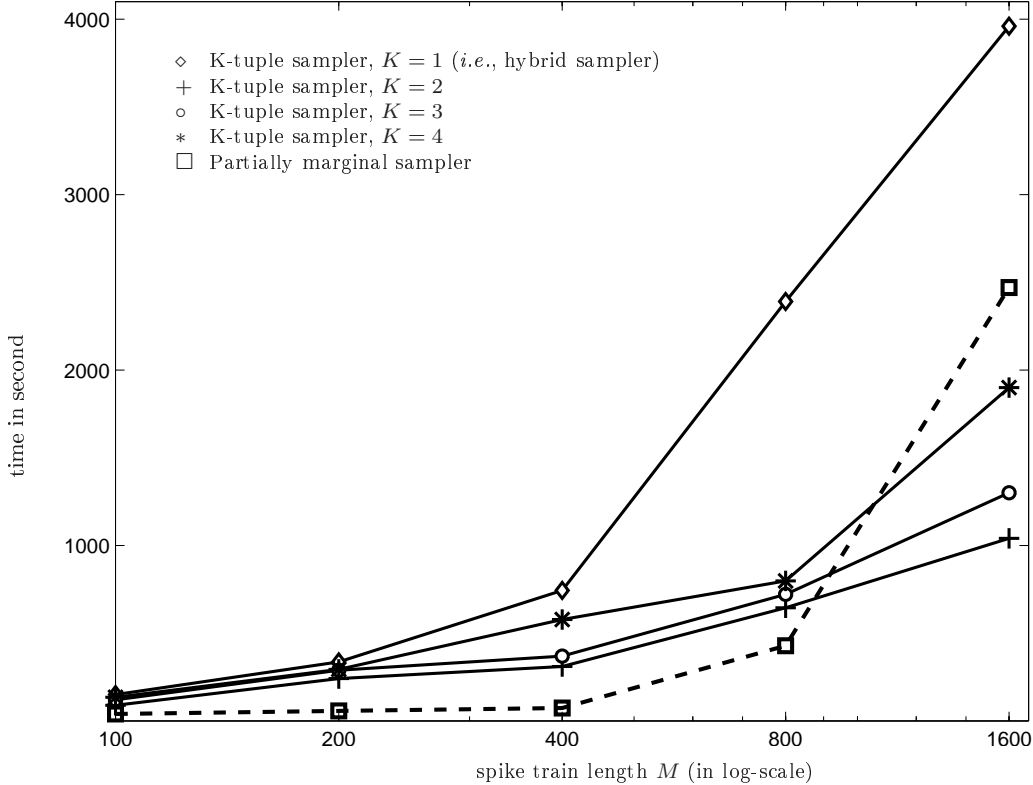}}%
{\hspace{60mm} {\footnotesize spike train length $M$ (in log-scale)}}
{\hspace{40mm}{\footnotesize time in second}}
{\hspace{20mm}\vspace{80mm}
	\renewcommand{\arraystretch}{0.8}
	\begin{tabular}{@{}l@{~}l@{}}
		{\scriptsize $\diamond$}        &{\scriptsize K-tuple sampler, $K=1$ (\ie hybrid sampler)}\\
		{\scriptsize $+$ }              &{\scriptsize K-tuple sampler, $K=2$}\\
		{\scriptsize $\circ$  }   			&{\scriptsize K-tuple sampler, $K=3$} \\
		{\scriptsize $\ast$}            &{\scriptsize K-tuple sampler, $K=4$} \\
		{\scriptsize $\square$}         &{\scriptsize Partially marginal sampler}
	\end{tabular}
}
\caption{Convergence time comparison for the $5$ samplers applied on spike trains of increasing length. \label{fig:ConvergTime}}
\end{figure}

In the given example of the Mendel's sequence (for which $M=300$), the partially marginal sampler achieves a better compromise between the two criteria~\cite{Liu01} than the most performant sampler in the generalized $K$-tuple family: the time required to reach convergence is reduced by half.
However, the partially marginal sampler do not necessarily outperform the $K$-tuple sampler in dealing with longer datasets.
Notice that the computational complexity for a Cholesky update operation in the latter case is $\Oc(L^2)$ ($L=\sum q_i\approx \lambda M$, the number of detected spikes), thus proportional to the data length of the BG sequence, whereas the complexity per iteration of the $K$-tuple sampler is non-increasing with respect to the number of detected spikes $L$. To illustrate the disadvantage of the partially marginal sampler, it is necessary to test on simulation examples of longer observation data \zb while fixing the Bernoulli parameter $\lambda$. Figure~\ref{fig:CostIter} compares the complexity per iteration as a function of BG sequence length $M$ for the two algorithms. The $2$-tuple sampler potentially outperforms the partially marginal sampler due to its linear growth of complexity per iteration with respect to $M$. 
Figure~\ref{fig:ConvergTime} resumes the convergence time for all the $5$ samplers ($K=1,\ldots, 4$ in solide lines and the partially marginalised sampler in dashed line) on simulation data whose lengths vary from $100$ to $1600$ and $\lambda = 10\%$. While the partially marginal sampler outperforms all the rest algorithms in areas of shorter data length ($M<800$), its quadratic computation cost per iteration as shown in Fig.~\ref{fig:CostIter} yields the contrary in applications with longer data lengths. We also note that for the given bernoulli parameter $\lambda = 10\%$, the hybrid sampler is never a good choice inside the $K$-tuple family regardless of the data length. 
  
\section{Conclusion}\label{sec:Conclusion}
This paper proposes two distinct methods in dealing with the identified inefficiency concerning the Bernoulli label \qb sampling in the BG deconvolution problem. Detailed algorithms are given for both methods as well as their performance on simulation tests. 
It is shown that both methods are mathematically valid MCMC sampling schemes and both achieve better convergence properties in comparison to the hybrid sampler, the latter is based on Cheng \etal{}'s Gibbs scheme with time-shift compensations and scale re-sampling. Both proposed methods demonstrate better trade-off for the convergence rate by reducing drastically the iteration steps needed to converge. On the simulation test of Mendel's sequence, the partially marginal sampler is shown to out-perform the whole class of K-tuple samplers, among which the optimized parameter is found at $K=2$. This is however due to the relatively small size of the simulation problem and the contrary is observed by simply augmenting the simulated data length. 
 
We conclude with some extensions of the proposed techniques. By combining the idea of a K-tuple joint sampling and that of partial marginalization of \xb, \ie updating $\qb_{i+\Omega} \I \qb_{-(i+\Omega)},\hb,\sigma_{\epsilon},\lambda$ in Step~1 of Table~\ref{Di}, one might further reduce the number of iterations necessary to escape from local optima at the cost of higher computation load per iteration. Such a solution should be adopted when the deconvolution problem presents more difficulties, characterized by either low SNR, band-limited IR and/or less sparse spike trains. 

\renewcommand{\theequation}{A-\arabic{equation}}    
\setcounter{equation}{0}  
\section*{APPENDIX}
To update the Cholesky factor \Rv on removing a spike at the location $i$ such that $\Kv\T\Kv = \Rv\T\Rv - \vb\vb\T$, we specify : 
\bal
\Rv = \left[
\barr{ccc}
\Av & \bb & \Cv \\
0   & \delta & \eb\T \\
0   &  0     & \Fv  
\earr
\right], 
\Kv = \left[
\barr{ccc}
\Gv & 0 & \Hv \\
0   & 0 & 0 \\
0   &  0     & \Jv  
\earr
\right],\label{eq:RK}\\
\left\{
\barr{ll}
\vb &= \sqrt{\tau}\Rv\T \cdot \Rv(:,j) \\
\tau\M &= \bb\T \bb + \delta^2
\earr
\right. \label{eq:vtau}
\eal
where $\{\Av, \Fv, \Gv, \Jv\}$ are upper-triangular matrices and $\delta = \Rv_{i,i}$ the $i$th diagonal value of \Rv. $(\Kv_{-i,-i})$ is therefor the updated Cholesky factor. 

A most evident solution consists an operation $\Kv = \text{Cholupdate}(\Rv, \vb, '-')$ followed by an extraction to get the updated Cholesky factor (of lowered dimension) :
\balx
\left[
\barr{cc}
\Gv & \Hv \\
0   & \Jv  
\earr
\right], 
\ealx
which however systematically fails due to the fact that $\Kv\T \Kv$ is no longer definite positive. In fact, one necessary condition~\cite{Seeger08} for a successful rank-$1$ \emph{downdate} operation (update by subtracting $\vb\vb\T$) is that $\rho^2 = 1- \pb\T\pb > 0$ for which \pb is defined by 
\bal
\Rv\T \pb = \vb. 
\eal
Without inverting $\Rv\T$, it can be deduced from ~\eqref{eq:vtau} that $\pb = \sqrt{\tau} \Rv(:,j)$ and thus $\rho^2 = 0$. The algorithm of \emph{Cholupdate} by a series of \emph{Givens} rotation should break down according to~\cite{Seeger08}. 

The problem of non-definite positiveness is avoided by a strategy that first lower the dimension and then perform the \emph{Cholupdate} algorithm. 
From ~\eqref{eq:RK}, we obtain: 
\balx
\Rv\T \Rv = \left[
\barr{c|c|c}
\Av\T\Av & \Av\T\bb & \Av\T\Cv \\ \hline
\bb\T\Av & \tau\M   & \bb\T\Cv + \delta\eb\T \\ \hline
\Cv\T\Av & \Cv\T\bb + \delta\eb   & \Cv\T\Cv + \eb\eb\T + \Fv\T\Fv  
\earr
\right]
= \left[
\barr{c|c|c}
\Gv\T\Gv & 0 & \Gv\T\Hv \\ \hline
0   & 0 & 0 \\ \hline
\Hv\T\Gv   &  0     & \Hv\T\Hv +\Jv\T\Jv   
\earr
\right]+ \vb\vb\T
\ealx
The submatrices on both sides by eliminating the $i$th row and column yields :
\balx
\left[
\barr{c|c}
\Av\T\Av &  \Av\T\Cv \\ \hline
\Cv\T\Av & \Cv\T\Cv + \eb\eb\T + \Fv\T\Fv  
\earr
\right]
= \left[
\barr{c|c}
\Gv\T\Gv  & \Gv\T\Hv \\  \hline
\Hv\T\Gv  & \Hv\T\Hv +\Jv\T\Jv   
\earr
\right]+ \vb_{-i}\vb_{-i}\T
\ealx
for which $\vb_{-i}$ denotes \vb without the $i$th component. 
The following relation is key to the algorithm that we propose : 
\bal
\left[
\barr{cc}
\Av &  \Cv \\
0   &  \Fv
\earr
\right]\T 
\left[
\barr{cc}
\Av &  \Cv \\
0   &  \Fv
\earr
\right] + \left[
\barr{c}
0 \\ \eb
\earr
\right] \left[
\barr{c}
0 \\ \eb
\earr
\right]\T
= \left[
\barr{cc}
\Gv  & \Hv \\
0    & \Jv   
\earr
\right]\T \left[
\barr{cc}
\Gv  & \Hv \\
0    & \Jv   
\earr
\right] + \vb_{-i}\vb_{-i}\T
\eal
The two steps that leads to a spike removing update are : 
\ben[align=left, leftmargin=*, noitemsep]
\item an update on \Fv (of reduced dimension): 
\balx
\Mv = \text{Cholupdate}(\Fv, \eb, '+');  
\ealx
\item a downdate that keeps the positiveness of the product matrix and will not fail: 
\balx
\left[
\barr{cc}
\Gv  & \Hv \\
0    & \Jv   
\earr
\right] = \text{Cholupdate}(\left[
\barr{cc}
\Av &  \Cv \\
0   &  \Mv
\earr
\right], \vb_{-i}, '-'). 
\ealx
it can be proven that $\rho^2 = 1-\pb\T\pb$ is strictly positive in this case. 
\een

\bibliographystyle{ieeeji}

\end{document}

%% file: jidef_en.tex
\def\dd{\,d}

\def\eqname{Eq.~}
\def\eqnames{Eqs~}
\def\tabname{Table~}

\def\figname{Figure~}

\def\chapname{Chapter~}

\def\parname{\S~}
\def\parnames{\S~}
\def\sectname{Section~}

\def\<{``}\def\>{''}

\def\etal{{\itshape et al.}}


%% file: beginend.tex
\def\barr{\begin{array}}                \def\earr{\end{array}}
\def\bcc{\begin{center}}                \def\ecc{\end{center}}
\def\cl#1{\centerline{#1}}
\def\bdoc{\begin{document}}             \def\edoc{\end{document}}
\def\ben{\begin{enumerate}}             \def\een{\end{enumerate}}
\def\beqn{\begin{eqnarray}}             \def\eeqn{\end{eqnarray}}
\def\beqnl#1{\beqn\label{#1}}           \def\eeqnl#1{\label{#1}\eeqn}
\def\beqnx{\begin{eqnarray*}}           \def\eeqnx{\end{eqnarray*}}
\def\bseqn{\begin{subeqnarray}}         \def\eseqn{\end{subeqnarray}}
\def\beq#1\eeq{\begin{equation}#1\end{equation}}
\def\bal#1\eal{\begin{align}#1\end{align}}
\def\balx#1\ealx{\begin{align*}#1\end{align*}}
\def\beqx{$$}                           \def\eeqx{$$}
\def\bfig{\protect\begin{figure}}       \def\efig{\protect\end{figure}}
\def\bfigx{\protect\begin{figure*}}     \def\efigx{\protect\end{figure*}}
\def\bfigt{\protect\begin{figurette}}   \def\efigt{\protect\end{figurette}}
\def\bfl{\begin{flushleft}}             \def\efl{\end{flushleft}}
\def\bfr{\begin{flushright}}            \def\efr{\end{flushright}}
\def\bit{\begin{itemize}}               \def\eit{\end{itemize}}
\def\bmi{\begin{minipage}}              \def\emi{\end{minipage}}
\newcommand{\bfmi}{\begin{fminipage}}   \newcommand{\efmi}{\end{fminipage}}
\def\bpic{\begin{picture}}              \def\epic{\end{picture}}
\def\btabl{\begin{table}}               \def\etabl{\end{table}}
\def\btab{\begin{tabular}} 
\def\btabu{\begin{tabular}}             \def\etabu{\end{tabular}}
\def\btabx{\begin{tabular*}}            \def\etabx{\end{tabular*}}

%% file: sigproc09.bbl
\begin{thebibliography}{}

\end{thebibliography}


\begin{thebibliography}{10}

\bibitem{Cheng96}
Q.~Cheng, R.~Chen, and T.-H. Li,
\newblock ``Simultaneous wavelet estimation and deconvolution of reflection
  seismic signals'',
\newblock {\em {{IEEE} {T}rans. {G}eosci. {R}emote {S}ensing}}, vol. 34, pp.
  377--384, {M}ar. 1996.

\bibitem{Kormylo83a}
J.~J. Kormylo and J.~M. Mendel,
\newblock ``Maximum-likelihood seismic deconvolution'',
\newblock {\em {{IEEE} {T}rans. {G}eosci. {R}emote {S}ensing}}, vol. GE-21, no.
  1, pp. 72--82, {J}an. 1983.

\bibitem{Bourguignon05}
S.~Bourguignon and H.~Carfantan,
\newblock ``Bernoulli-{G}aussian spectral analysis of unevenly spaced
  astrophysical data'',
\newblock in {\em IEEE Workshop Stat. Sig. Proc.}, Bordeaux, France, {J}uly
  2005, pp. 811--816.

\bibitem{Champagnat96a}
F.~Champagnat, Y.~Goussard, and J.~Idier,
\newblock ``Unsupervised deconvolution of sparse spike trains using stochastic
  approximation'',
\newblock {\em {{IEEE} {T}rans. {S}ignal {P}rocessing}}, vol. 44, no. 12, pp.
  2988--2998, {D}ec. 1996.

\bibitem{Liu01}
J.~S. Liu,
\newblock {\em Monte {C}arlo Strategies in Scientific Computing},
\newblock Springer Series in Statistics. Springer Verlag, New York, NY, 2001.

\bibitem{Robert04}
C.~P. Robert and G.~Casella,
\newblock {\em Monte {C}arlo Statistical Methods},
\newblock Springer Texts in Statistics. Springer Verlag, New York, NY, 2nd
  edition, 2004.

\bibitem{Labat06b}
C.~Labat and J.~Idier,
\newblock ``Sparse blind deconvolution accounting for time-shift ambiguity'',
\newblock in {\em {{P}roc. {IEEE} {ICASSP}}}, Toulouse, France, {M}ay 2006, pp.
  616--619.

\bibitem{Veit08}
T.~Veit, J.~Idier, and S.~Moussaoui,
\newblock ``{R}ééchantillonnage de l'échelle dans les algorithmes {MCMC} pour
  les problèmes inverses bilinéaires'',
\newblock {\em {{T}raitement du {S}ignal}}, vol. 25, no. 4, pp. 329--343, 2008.

\bibitem{Chi84}
C.~Y. Chi and J.~M. Mendel,
\newblock ``Improved maximum-likelihood detection and estimation of
  {B}ernoulli-{G}aussian processes'',
\newblock {\em {{IEEE} {T}rans. {I}nf. {T}heory}}, vol. 30, pp. 429--435, 1984.

\bibitem{Ge08}
D.~Ge, J.~Idier, and E.~{Le Carpentier},
\newblock ``A new {MCMC} algorithm for blind {B}ernoulli-{G}aussian
  deconvolution'',
\newblock in {\em EUSIPCO}, Lausanne, Switzerland, {S}ep. 2008.

\bibitem{Dobigeon07}
N.~Dobigeon, J.-Y. Tourneret, and J.~D. Scargle,
\newblock ``Joint segmentation of multivariate astronomical times series:
  {B}ayesian sampling with a hierarchical model'',
\newblock {\em {{IEEE} {T}rans. {S}ignal {P}rocessing}}, vol. 55, no. 2, pp.
  414--423, {F}eb. 2007.

\bibitem{Brooks98}
S.~P. Brooks and A.~Gelman,
\newblock ``General methods for monitoring convergence of iterative
  simulations'',
\newblock {\em {{J}. {C}omput. {G}raph. {S}tatist.}}, vol. 7, no. 4, pp.
  434--455, 1998.

\bibitem{Liu94b}
J.~S. Liu, W.~H. Wong, and A.~Kong,
\newblock ``Covariance structure of the {G}ibbs sampler with applications to
  the comparisons of estimators and augmentation schemes'',
\newblock {\em {{B}iometrika}}, vol. 81, pp. 27--40, 1994.

\bibitem{Cappe99}
O.~Cappé, A.~Doucet, M.~Lavielle, and E.~Moulines,
\newblock ``Simulation-based methods for blind maximum-likelihood filter
  identification'',
\newblock {\em {{S}ignal {P}rocessing}}, vol. 73, no. 1-2, pp. 3--25, {F}eb.
  1999.

\bibitem{Park09}
T.~Park and D.~A. van Dyk,
\newblock ``{Partially Collapsed Gibbs Samplers: Illustrations and
  Applications}'',
\newblock {\em Journal of Computational and Graphical Statistics}, vol. 18, no.
  2, pp. 283--305, 2009.

\bibitem{Seeger08}
M.~Seeger,
\newblock ``{Low Rank Updates for the Cholesky Decomposition}'',
\newblock Tech. {R}ep., Department of EECS, University of California at
  Berkeley, 485 Soda Hall, Berkeley CA, 2008.

\end{thebibliography}
